\documentclass[twoside,11pt]{article}

\usepackage{jmlr2e}
\RequirePackage{amsmath}

\usepackage{lastpage}
\jmlrheading{22}{2021}{1-\pageref{LastPage}}{9/20; Revised
1/21}{3/21}{20-1074}{Lasse Petersen and Niels Richard Hansen}
\ShortHeadings{Testing Conditional Independence via Quantile Regression}{Petersen and Hansen}

\newtheorem{myassump}{Assumption}

\renewcommand{\d}{\mathrm{d}}

\usepackage{centernot}
\newcommand{\indep}{\mathrel{\text{\scalebox{1.07}{$\perp\mkern-10mu\perp$}}}}
\newcommand{\notindep}{\centernot{\indep}}

\DeclareMathOperator*{\argmin}{arg\,min}

\newcommand{\R}{\mathbb{R}}

\renewcommand{\epsilon}{\varepsilon}
\renewcommand{\phi}{\varphi}

\firstpageno{1}

\begin{document}

\title{Testing Conditional Independence via \\
Quantile Regression Based Partial Copulas}

\author{\name Lasse Petersen \email lp@math.ku.dk \\
       \addr Department of Mathematical Sciences\\
       University of Copenhagen\\
       Universitetsparken 5, 2100 Copenhagen, Denmark
       \AND
       \name Niels Richard Hansen \email niels.r.hansen@math.ku.dk \\
       \addr Department of Mathematical Sciences\\
       University of Copenhagen\\
       Universitetsparken 5, 2100 Copenhagen, Denmark}

\editor{Peter Spirtes}

\maketitle

\begin{abstract}
The partial copula provides a method for describing the dependence between two
random variables $X$ and $Y$ conditional on a third random vector $Z$ in 
terms of nonparametric residuals $U_1$ and $U_2$. This paper develops a nonparametric 
test for conditional independence by combining the partial copula with a  
quantile regression based method for estimating the nonparametric residuals.
We consider a test statistic based on generalized correlation between $U_1$ and $U_2$ and 
derive its large sample properties under consistency assumptions on the quantile regression 
procedure. We demonstrate through a simulation study 
that the resulting test is sound under complicated 
data generating distributions. 
Moreover, in the examples considered the test is competitive to other state-of-the-art 
conditional independence tests in terms of level and power, and it has superior power 
in cases with conditional variance heterogeneity of $X$ and $Y$ given $Z$.
\end{abstract}


\section{Introduction} \label{sec:introduction}

This paper introduces a new class of nonparametric tests of conditional 
independence between real-valued random variables, $X \indep Y \mid Z$, based on quantile regression. Conditional independence is an important concept in many statistical 
fields such an graphical models and causal inference 
\citep{lauritzen1996, spirtes2000, pearl2009}. However, 
\cite{shah2018hardness} proved that conditional independence is an untestable 
hypothesis when the distribution of $(X, Y, Z)$ is only assumed to be 
absolutely continuous with respect to Lebesgue measure.

More precisely, let $\mathcal{P}$ denote the set of distributions of 
$(X, Y, Z)$ that are absolutely continuous 
with respect to Lebesgue measure. Let $\mathcal{H} \subset \mathcal{P}$ be those
distributions for which conditional independence holds. 
Then \cite{shah2018hardness} showed that if $\psi_n$ is a hypothesis test 
for conditional independence with uniformly 
valid level $\alpha \in (0, 1)$ over $\mathcal{H}$,
  \begin{align*}
  \sup_{P \in \mathcal{H}} E_P(\psi_n) \leq \alpha,
  \end{align*}
then the test cannot have power greater than 
$\alpha$ against any alternative $P \in \mathcal{Q} := \mathcal{P} \setminus 
\mathcal{H}$. This is true even when restricting the distribution of $(X, Y, Z)$ to have
bounded support. The purpose of 
this paper is to identify a subset $\mathcal{P}_0 \subset \mathcal{P}$ of distributions 
and a test $\psi_n$ that has asymptotic (uniform)
level over $\mathcal{P}_0 \cap \mathcal{H}$ and power 
against a large set of alternatives within $\mathcal{P}_0 \backslash \mathcal{H}$.

Our starting point is the 
so-called partial copula construction. 
Letting $F_{X \mid Z}$ and $F_{Y \mid Z}$ denote
the conditional distribution functions of $X$ given $Z$ and $Y$ given $Z$, 
respectively, we define random variables $U_1$ and $U_2$ by
  \begin{align*}
  U_1 := F_{X \mid Z}(X \mid Z) \quad \text{and} \quad 
  U_2 := F_{Y \mid Z}(Y \mid Z).
  \end{align*}
Then the joint distribution of $U_1$ and $U_2$ is called the partial copula and
it can be shown that $X \indep Y \mid Z$ implies $U_1 \indep U_2$. Thus the 
question about conditional independence can be transformed into a question 
about independence. The main challenge with this approach is that 
the conditional distribution functions are unknown and 
must be estimated. 

In Section \ref{sec:conddistest} we
propose an estimator of conditional 
distribution functions based on quantile regression. 
More specifically, we let $\mathcal{T} = [\tau_{\min}, \tau_{\max}]$ 
be a range of quantile levels for 
$0 < \tau_{\min} < \tau_{\max} < 1$, and let 
$Q(\mathcal{T} \mid z)$ denote the range of conditional 
$\mathcal{T}$-quantiles in the distribution $X \mid Z=z$. 
To estimate a conditional distribution 
function $F$ given a sample $(X_i, Z_i)_{i=1}^n$ we propose 
to perform quantile regressions $\hat{q}_{k, z} = 
\hat{Q}^{(n)}(\tau_k \mid z)$ along an equidistant grid of quantile levels 
$(\tau_k)_{k=1}^m$ in $\mathcal{T}$, and then construct the estimator 
$\hat{F}^{(m, n)}$ by linear interpolation of the points 
\mbox{$(\hat{q}_{k, z}, \tau_k)_{k=1}^m$}. The main result of 
the first part of the paper is Theorem \ref{thm:pointwiseconsistency}, 
which states that we can achieve the following bound on the estimation 
error
	\begin{align*}
	\| F - \hat{F}^{(m, n)} \|_{\mathcal{T}, \infty} :=
	\sup_{z} \sup_{t \in Q(\mathcal{T} \mid z)}
	|F(t \mid z) - \hat{F}^{(m, n)}(t \mid z) |
	\in \mathcal{O}_P(g_P(n)) 
	\end{align*}
where $g_P$ is a rate function describing the 
$\mathcal{O}_P$-consistency of the quantile regression procedure 
over the conditional $\mathcal{T}$-quantiles 
for $P$ in a specified set of distributions $\mathcal{P}_0 \subset 
\mathcal{P}$. This result demonstrates 
how pointwise consistency of a quantile regression procedure 
over $\mathcal{P}_0$ can be transferred 
to the estimator $\hat{F}^{(m, n)}$, and we discuss how this 
can be extended to uniform consistency over $\mathcal{P}_0$. 
We conclude the section by 
reviewing a flexible model class from quantile regression where 
such consistency results are available.

In Section \ref{sec:test-cond-ind} we describe a generic method 
for testing conditional independence based on  
estimated conditional distribution functions, 
$\hat{F}_{X \mid Z}^{(n)}$ and $\hat{F}^{(n)}_{Y \mid Z}$, 
obtained from a sample $(X_i, Y_i, Z_i)_{i=1}^n$.
From these estimates we compute 
  \begin{align*}
  \hat{U}_{1, i}^{(n)} := \hat{F}^{(n)}_{X \mid Z}(X_i \mid Z_i) \quad \text{and} 
  \quad \hat{U}_{2, i}^{(n)} := \hat{F}^{(n)}_{Y \mid Z}(Y_i \mid Z_i),
  \end{align*}
for $i=1, \dots, n$, which can then be plugged into a
bivariate independence test. If \mbox{$\hat{F}_{X \mid Z}^{(n)}$} and
\mbox{$\hat{F}^{(n)}_{Y \mid Z}$} are consistent with a sufficiently fast 
rate of convergence, properties of the bivariate test, in terms of
level and power, can be transferred to the test of conditional
independence. 
The details of this transfer of properties depend on the specific test
statistic. The main contribution of the second part of the paper is a
detailed treatment of a test given in terms of a generalized 
correlation, estimated as 
  \begin{align*}
  \hat{\rho}_n := \frac{1}{n} \sum_{i = 1}^n 
  \phi \left( \hat{U}_{1, i}^{(n)} \right) 
  \phi \left( \hat{U}_{2, i}^{(n)} \right)^T
  \end{align*}
for a function 
$\phi = (\phi_1, \dots, \phi_q): [0,1] \to \mathbb{R}^q$ satisfying 
certain regularity conditions. A main result is Theorem 
\ref{thm:rhoasympnormal}, which states that $\sqrt{n} \hat{\rho}_n$ 
converges in distribution toward $\mathcal{N}(0, \Sigma \otimes \Sigma)$ under the  
hypothesis of conditional independence whenever 
$\hat{F}_{X \mid Z}^{(n)}$ and $\hat{F}^{(n)}_{Y \mid Z}$ are $\mathcal{O}_P$-consistent with 
rates $g_P$ and $h_P$ satisfying $\sqrt{n} g_P(n) h_P(n) \to 0$. 
The covariance matrix $\Sigma$ depends only on $\phi$.
We use this 
to show asymptotic pointwise level of the test when restricting to the set of distributions 
$\mathcal{P}_0$ where the required consistency can be obtained.
We then proceed to show in Theorem \ref{thm:rhohat_to_rho} that $\sqrt{n} \hat{\rho}_n$ 
diverges in probability under a set of alternatives of 
conditional dependence when we have $\mathcal{O}_P$-consistency of the 
conditional distribution function estimators.
This we use to show asymptotic pointwise power of the test. 
We also show how asymptotic uniform level and power can be achieved
when $\hat{F}_{X \mid Z}^{(n)}$ and $\hat{F}^{(n)}_{Y \mid Z}$  
are uniformly consistent over $\mathcal{P}_0$. Lastly, we 
provide an out-of-the-box procedure for 
conditional independence testing in conjunction with our 
quantile regression based conditional distribution function 
estimator $\hat{F}^{(m, n)}$ from Section \ref{sec:conddistest}.

In Section \ref{sec:simulations} we examine the proposed test 
through a simulation study 
where we assess the level and power properties of the test and 
benchmark it against 
existing nonparametric conditional independence tests. 
All proofs are collected in Appendix \ref{appendix:proofs}.

\section{Related Work} \label{sec:related_work}

The partial copula and its application for conditional independence testing was initially 
introduced by \cite{bergsma2004testing} and further explored by \cite{bergsma_nonparametric_2011}. 
Its use for conditional independence testing has also been explored by \cite{song_testing_2009}, \cite{patra2016nonparametric} 
and \cite{liu2018}. Moreover, properties of the partial copula was 
studied by \cite{gijbels_estimation_2015} and \cite{spanhel_partial_2016} among others. 
A related but different approach for testing conditional independence via the factorization
of the joint copula of $(X, Y, Z)$ is given by 
\cite{bouezmarni2012nonparametric}. 
Common for the existing approaches to using the partial copula for conditional independence testing 
is that the conditional distribution functions $F_{X \mid Z}$ and $F_{Y \mid Z}$ are estimated 
using a kernel smoothing procedure \citep{stute1986almost, einmahl2005uniform}. 
The advantage of the approach is that the estimator is nonparametric, however, 
it does not scale well with the dimension of the conditioning variable $Z$. 
This is partly remedied by \cite{haff2015nonparametric} who suggest a nonparametric pair-copula 
estimator whose convergence rate is independent of the dimension of $Z$. This estimator 
requires the simplifying assumption, which is a strong assumption not 
required for the validity of our approach. Moreover, it is not obvious how to incorporate 
parametric assumptions, such as a certain functional dependence between 
response and covariates, using kernel smoothing estimators, since there is only the choice of a 
kernel and a bandwidth.
Furthermore, a treatment of the relationship between level and power properties 
of a partial copula based conditional independence test, and consistency 
of the conditional distribution function estimator is lacking in the existing literature. 
In this work we take a novel approach to testing conditional independence using the partial copula 
by using quantile regression for estimating the conditional distribution functions. This allows for 
a distribution free modeling of the conditional distributions $X \mid Z=z$ and $Y \mid Z=z$ that 
can handle high-dimensionality of $Z$ through penalization, and complicated response-predictor 
relationships by basis expansions. We also make the requirements on consistency 
of the conditional distribution function estimator 
that are needed to obtain level and power of the test explicit. A similar recent
approach to testing conditional independence using regression methods is given by \cite{shah2018hardness}, who 
propose to test for vanishing correlation between the residuals after nonparametric conditional 
mean regression of $X$ on $Z$ and $Y$ on $Z$. See also \cite{ramsey2014scalable} and 
\cite{fan2020projection}. This approach captures dependence between 
$X$ and $Y$ given $Z$ that lies in the conditional correlation. However, as is demonstrated through a simulation study in 
Section \ref{sec:local_alternative}, it does not adequately 
account for conditional variance heterogeneity between $X$ and $Y$ given $Z$, 
while our partial copula based test captures the dependence more efficiently.

\section{Estimation of Conditional Distribution Functions} \label{sec:conddistest}

Throughout the paper we restrict ourselves to the set of distributions $\mathcal{P}$ over the 
hypercube $[0, 1]^{2+d}$ that are absolutely continuous with respect to 
Lebesgue measure. Let $(X, Y, Z) \sim P \in \mathcal{P}$ such that 
$X, Y \in [0, 1]$ and $Z \in [0, 1]^d$. When we speak of the distribution 
of $X$ given $Z$ relative to $P$ we mean the conditional distribution 
that is induced when $(X, Y, Z) \sim P$. 
In this section we consider 
estimation of the conditional distribution function $F_{X \mid Z}$ of 
$X$ given $Z$ using quantile regression. Estimation of $F_{Y \mid Z}$ can be carried out analogously. 

\subsection{Conditional distribution and quantile functions}

Given $z \in [0, 1]^d$ we denote by
	\begin{align*}
	F_{X \mid Z}(t \mid z) := P(X \leq t \mid Z=z)
	\end{align*}
the conditional distribution function of $X \mid Z=z$ for $t \in [0, 1]$. 
We denote by 
	\begin{align*}
	Q_{X \mid Z} (\tau \mid z) := \inf\{ t \in [0, 1] 
	\mid F_{X\mid Z}(t \mid z) \geq \tau\}
	\end{align*}
the conditional quantile function of the conditional distribution 
$X \mid Z=z$ for $\tau \in [0, 1]$ and $z \in [0, 1]^d$. 
We will omit the subscript in $F_{X \mid Z}$ and $Q_{X \mid Z}$
when the conditional distribution of interest is clear from the context. 

In quantile regression one models the function 
$z \mapsto Q(\tau \mid z)$ for fixed $\tau \in [0, 1]$.
Estimation of the quantile regression function
is carried out by solving the empirical risk minimization problem
  \begin{align*}
  \hat{Q}(\tau \mid \cdot) \in \argmin_{f \in \mathcal{F}} \sum_{i=1}^n
  L_\tau\left( 
  X_i - f(Z_i)
  \right) 
  \end{align*}
where the loss function $L_\tau(u) =  u(\tau - 1(u < 0))$ is the 
so-called check function and $\mathcal{F}$ is some function class. 
For $\tau = 1 / 2$ the loss function is $L_{1 / 2}(u) = |u|$, and 
we recover median regression as a special case. One can also choose to 
add regularization as with conditional mean regression.
See \cite{quantileregression} and \cite{koenker2017handbook} for 
an overview of the field.

\subsection{Quantile regression based estimator} \label{sec:quantregbasedestimator}

Based on the conditional quantile function $Q$ we define an
approximation $\tilde{F}^{(m)}$ of the conditional distribution function $F$ as
follows. We let $\tau_{\min}$ and $\tau_{\max}$ denote fixed quantile levels 
satisfying $0 < \tau_{\min} < \tau_{\max} < 1$, and we let 
$q_{\min, z} := Q(\tau_{\min} \mid z) > 0$ and $q_{\max, z} := 
Q(\tau_{\max} \mid z) < 1$ denote the corresponding conditional quantiles. 

Let $\mathcal{T} = [\tau_{\min}, \tau_{\max}]$ 
denote the set of potential quantile levels. 
A grid in $\mathcal{T}$ is a sequence $(\tau_k)_{k=1}^m$ such that 
$\tau_{\min} = \tau_1 < \cdots < \tau_{m} = \tau_{\max}$ for $m \geq
2$. An equidistant grid is a grid $(\tau_k)_{k=1}^m$
for which $\tau_{k+1} - \tau_k$ is constant for $k=1, \dots, m-1$. 
Also let $\tau_0 = 0$ and $\tau_{m+1} = 1$ be fixed. 

Given a grid $(\tau_k)_{k=1}^m$ we let \mbox{$q_{k, z} := Q(\tau_k \mid z)$} for 
$k = 1, \dots, m$ and define $q_{0, z} := 0$ and $q_{m+1, z} := 1$. 
For each $z \in [0, 1]^d$ we define a function 
$\tilde{F}^{(m)}(\cdot \mid z) : [0, 1] \to [0, 1]$ 
by linear interpolation of the points $(q_{k, z}, \tau_k)_{k=0}^{m+1}$:
  \begin{align} \label{eq:Ftilde}
  \tilde{F}^{(m)}(t \mid z) 
  :=
  \sum_{k=0}^{m}
  \left(
  \tau_k + (\tau_{k+1} - \tau_k) \frac{t - q_{k, z}}{q_{k+1, z} - q_{k, z}}
  \right)
  1_{(q_{k, z}, q_{k+1, z}]}(t).
  \end{align}
Let $Q(\mathcal{T} \mid z) = [q_{\min, z}, q_{\max, z}]$ be the 
range of conditional $\mathcal{T}$-quantiles in the conditional 
distribution $X \mid Z=z$ for $z \in [0, 1]^d$, and define  
the supremum norm
	\begin{align*}
	\| f \|_{\mathcal{T}, \infty} = 
	\sup_{z \in [0, 1]^d} \sup_{t \in Q(\mathcal{T} \mid z)} 
	|f(t, z)|
	\end{align*}
for a function $f : [0, 1] \times [0, 1]^d \to \R$. Note that 
this is a norm on the set of bounded functions on 
$\{ (t, z) \mid z \in [0,1]^d, t \in Q(\mathcal{T} \mid z) \}$.
Then we have the following approximation result.

\begin{proposition} \label{prop:Ftildeapproximation}
Denote by $\tilde{F}^{(m)}$ the function \eqref{eq:Ftilde} defined 
from a grid $(\tau_k)_{k=1}^m$ in $\mathcal{T}$. Then it holds 
that
	\begin{align*}
	||F - \tilde{F}^{(m)}||_{\mathcal{T}, \infty} \leq \kappa_m 
	\end{align*}
where $\kappa_m := \max_{k=1, \dots, m-1}(\tau_{k+1} - \tau_k)$ is 
the coarseness of the grid.
\end{proposition}

Choosing a finer and finer grid yields $\kappa_m \to 0$, which implies that 
$\tilde{F}^{(m)} \to F$ in the norm $\| \cdot \|_{\mathcal{T}, \infty}$ 
for $m \to \infty$.

By an estimator of the conditional distribution function $F$ we 
mean a mapping from a sample $(X_i, Z_i)_{i=1}^n$
to a function \mbox{$\hat{F}^{(n)}(\cdot \mid z) : [0, 1] \to [0, 1]$} 
such that for every 
$z \in [0, 1]^d$ it holds that 
$t \mapsto \hat{F}^{(n)}(t \mid z)$ is continuous and increasing with
  \begin{align*}
  \hat{F}^{(n)}(0 \mid z) = 0 \quad \text{and} \quad
  \hat{F}^{(n)}(1 \mid z) = 1.
  \end{align*}
Motivated by \eqref{eq:Ftilde} we define the following
estimator of the conditional distribution function.

\begin{definition} \label{def:Fhat}
Let $(\tau_k)_{k=1}^m$ be a grid in $\mathcal{T}$. Define $\hat{q}^{(n)}_{0, z} := 0$ and 
$\hat{q}_{m+1, z}^{(n)} := 1$, and let $\hat{q}^{(n)}_{k, z} := \hat{Q}^{(n)}(\tau_k \mid z)$
for $k=1, \dots, m$ be the predictions of a quantile regression model 
obtained from an i.i.d. sample $(X_i, Z_i)_{i=1}^n$.
We define the estimator 
$\hat{F}^{(m, n)}(\cdot \mid z) : [0, 1] \to [0, 1]$ by
  \begin{align} \label{eq:Fhat}
  \hat{F}^{(m, n)}(t \mid z) 
  :=
  \sum_{k=0}^{m}
  \left(
  \tau_k + (\tau_{k+1} - \tau_k)
  \frac{t - \hat{q}^{(n)}_{k, z}}{\hat{q}^{(n)}_{k+1, z} - \hat{q}^{(n)}_{k, z}}
  \right)
  1_{\left(\hat{q}^{(n)}_{k, z}, \hat{q}^{(n)}_{k+1, z} \right]}(t)
  \end{align}
for each $z \in [0, 1]^d$.
\end{definition}

Note that the estimator is not monotone in the presence of 
quantile crossing \citep{he1997quantile}. In this case we perform a 
re-arrangement of the estimated conditional quantiles in order to obtain monotonicity 
for finite sample size \citep{chernozhukov2010quantile}.
However, the estimated conditional quantiles will be ordered correctly under 
the consistency assumptions that we will introduce in Assumption 
\ref{assump:Qhatconsistency}, that is, the re-arrangement becomes unnecessary, 
and the estimator becomes monotone with high probability as $n \to \infty$ 
for any grid $(\tau_k)_{k=1}^m$ in $\mathcal{T}$.

\subsection{Pointwise consistency of $\hat{F}^{(m, n)}$} \label{sec:pointwiseconsistency}

We will now demonstrate how pointwise consistency of the proposed 
estimator over a set of distributions $\mathcal{P}_0 \subset \mathcal{P}$ 
can be obtained under the assumption that the quantile regression 
procedure is pointwise consistent over $\mathcal{P}_0$.

We will evaluate the consistency of $\hat{F}^{(m, n)}$ according to 
the supremum norm $||\cdot||_{\mathcal{T}, \infty}$ introduced in Section \ref{sec:quantregbasedestimator}, that is, we restrict the supremum to  
be over $t \in Q(\mathcal{T} \mid z)$ 
and not the entire interval $[0, 1]$. We do so because quantile regression
generally does not give uniform consistency of all extreme quantiles,
and in Section \ref{sec:test-cond-ind} we show how 
consistency of $\hat{F}^{(m, n)}$ between the conditional 
$\tau_{\min}$- and $\tau_{\max}$-quantiles is sufficient for conditional
independence testing. 

First, we have the following key 
corollary of Proposition \ref{prop:Ftildeapproximation}, which is a simple 
application of the triangle inequality.

\begin{corollary} \label{cor:keyinequality}
Let $\tilde{F}^{(m)}$ and $\hat{F}^{(m, n)}$ be given by
\eqref{eq:Ftilde} and \eqref{eq:Fhat}, respectively. Then
  \begin{align*}
  \|F -  \hat{F}^{(m, n)}\|_{\mathcal{T}, \infty}
  \leq 
  \kappa_m
  +
  \|\tilde{F}^{(m)} - \hat{F}^{(m, n)}\|_{\mathcal{T}, \infty}
  \end{align*}
for all grids $(\tau_k)_{k=1}^m$ in $\mathcal{T}$. 
\end{corollary}

The random part of the right hand side of the inequality is the term 
$\|\tilde{F}^{(m)} - \hat{F}^{(m, n)}\|_{\mathcal{T}, \infty}$, while $\kappa_m$ is deterministic 
and only depends on the choice of grid $(\tau_{k})_{k=1}^m$. Controlling the term 
$\|\tilde{F}^{(m)} - \hat{F}^{(m, n)}\|_{\mathcal{T}, \infty}$ is an easier task than controlling 
$\|F -  \hat{F}^{(m, n)}\|_{\mathcal{T}, \infty}$ directly because 
$\tilde{F}^{(m)}$ and $\hat{F}^{(m, n)}$ are piecewise linear, 
while $F$ is only assumed to be continuous and increasing.

Consistency assumptions on the quantile regression procedure will
allow us to show consistency of the estimator $\hat{F}^{(m, n)}$. Let the random variable
  \begin{align*}
  \mathcal{D}_{\mathcal{T}}^{(n)} :=
  \sup_{z \in [0, 1]^d}  \sup_{\tau \in \mathcal{T}}
  | Q(\tau \mid z) - \hat{Q}^{(n)}(\tau \mid z) |
  \end{align*}
denote the uniform prediction error of a fitted 
quantile regression model,  $\hat{Q}^{(n)}$, over the set of quantile levels
$\mathcal{T} = [\tau_{\min}, \tau_{\max}]$. Below we write 
$X_n \in \mathcal{O}_P(a_n)$ when $X_n$ is big-O in probability of 
$a_n$ with respect to $P$. See Appendix \ref{appendix:modesofconvergence} 
for the formal definition.

\pagebreak

\begin{myassump} \label{assump:Qhatconsistency} For each $P \in \mathcal{P}_0$ there exist
\begin{itemize}
\item[(i)] a deterministic rate function $g_P$ tending to zero as $n
  \to \infty$ such that $\mathcal{D}_{\mathcal{T}}^{(n)} \in \mathcal{O}_P(g_P(n))$
\item[(ii)] and a finite constant $C_P$ such that the conditional density $f_{X \mid Z}$ satisfies 
$$\sup_{x \in [0, 1]} f_{X \mid Z}(x \mid z) \leq C_P$$ for almost all 
$z \in [0, 1]^d.$
\end{itemize}
\end{myassump}

Assumption \ref{assump:Qhatconsistency} (i) is clearly necessary to 
achieve consistency of the estimator. Assumption
\ref{assump:Qhatconsistency} (ii) is 
a regularity condition that
is used to ensure that $q_{k+1, z} - q_{k, z}$ does not tend to zero too fast 
as $\kappa_m \to 0$. We now have:

\begin{proposition} \label{prop:Fhatconsistent} 
Let Assumption \ref{assump:Qhatconsistency} be satisfied. Then
  \begin{align*}
  \|\tilde{F}^{(m)} - \hat{F}^{(m, n)}\|_{\mathcal{T}, \infty} \in 
  \mathcal{O}_P(g_P(n))
  \end{align*}
for each fixed $P \in \mathcal{P}_0$ and all equidistant grids $(\tau_k)_{k = 1}^{m}$
in $\mathcal{T}$.
\end{proposition}

Consider letting the number of grid points $m_n$ depend on the sample size $n$. 
By combining Corollary \ref{cor:keyinequality} and Proposition 
\ref{prop:Fhatconsistent} we obtain the main pointwise consistency result. 

\begin{theorem} \label{thm:pointwiseconsistency}
Let Assumption \ref{assump:Qhatconsistency} be satisfied. Then 
  \begin{align*}
  \|F - \hat{F}^{(m_n, n)}\|_{\mathcal{T}, \infty} \in \mathcal{O}_P(g_P(n))
  \end{align*}
for each fixed $P \in \mathcal{P}_0$ given that the equidistant grids 
$(\tau_k)_{k = 1}^{m_n}$ in $\mathcal{T}$ satisfy $\kappa_{m_n} \in o(g_P(n))$.
\end{theorem}

This shows that $\hat{F}^{(m_n, n)}$ 
is pointwise consistent over $\mathcal{P}_0$ 
given that the 
quantile regression procedure is pointwise consistent over $\mathcal{P}_0$. 
Moreover, we can transfer the rate of convergence $g_P$ directly. 
In Section \ref{sec:correlationtest} we will use this type 
of pointwise consistency to 
show asymptotic pointwise level and power of our conditional independence 
test over $\mathcal{P}_0$.

Note that we can estimate conditional distribution functions in settings 
with high dimensional covariates to the extend that the quantile regression 
estimation procedure can deal with high dimensionality. An example of such a 
procedure is given in Section \ref{sec:quantileregmodels}.

We chose to state Proposition \ref{prop:Fhatconsistent} and 
Theorem \ref{thm:pointwiseconsistency} for equidistant grids only, but 
in the proof of Proposition \ref{prop:Fhatconsistent} we only need that 
the ratio $\kappa_m / \gamma_m$ between 
the coarseness $\kappa_m$ and the smallest 
subinterval $\gamma_m = \min_{k=1, \dots, m-1}(\tau_{k+1} - \tau_k)$ 
must not diverge as $m \to \infty$. This is obviously ensured for an equidistant 
grid. Moreover, for an equidistant grid,
$\kappa_m = (\tau_{\max} - \tau_{\min}) / (m-1)$, and $\kappa_{m_n} \in o(g_P(n))$
if $m_n$ grows with rate at least $g_P(n)^{-(1 + \epsilon)}$ for some
$\epsilon > 0$. Since the rate is unknown in practical applications 
we choose $m$ to be the smallest integer larger than $\sqrt{n}$ as a rule of thumb, 
since this represents the optimal parametric rate.

\subsection{Uniform consistency of $\hat{F}^{(m, n)}$}

The pointwise consistency result of Theorem \ref{thm:pointwiseconsistency} 
can be extended to a uniform consistency over $\mathcal{P}_0$ by strengthening 
Assumption \ref{assump:Qhatconsistency} to hold uniformly. 
Below we write $X_n \in \mathcal{O}_{\mathcal{M}}(a_n)$ when $X_n$ 
is big-O in probability of $a_n$ uniformly over a set of distributions 
$\mathcal{M}$. We refer to Appendix \ref{appendix:modesofconvergence} for the 
formal definition. 

\begin{myassump} \label{assump:uniformQhatconsistent} For
  $\mathcal{P}_0 \subset \mathcal{P}$ there exist
\begin{itemize}
\item[(i)] a deterministic rate function $g$ tending to 
zero as $n \to \infty$ such that $\mathcal{D}_{\mathcal{T}}^{(n)} \in \mathcal{O}_{\mathcal{P}_0}(g(n))$
\item[(ii)] and a finite constant $C$ such that the conditional density $f_{X \mid Z}$ satisfies 
$$\sup_{x \in [0, 1]} f_{X \mid Z}(x \mid z) \leq C$$ for almost all 
$z \in [0, 1]^d$.
\end{itemize}
\end{myassump}

With this stronger assumption we have a uniform 
extension of Proposition \ref{prop:Fhatconsistent}.

\begin{proposition} \label{prop:uniformFhatconsistent} 
Let Assumption \ref{assump:uniformQhatconsistent}
be satisfied. Then
  \begin{align*}
  \|\tilde{F}^{(m)} - \hat{F}^{(m, n)}\|_{\mathcal{T}, \infty} \in 
  \mathcal{O}_{\mathcal{P}_0}(g(n)).
  \end{align*}
for all equidistant grids $(\tau_k)_{k=1}^m$ in $\mathcal{T}$. 
\end{proposition}

We can now combine Corollary \ref{cor:keyinequality} 
with the stronger Proposition \ref{prop:uniformFhatconsistent}
to obtain the following uniform consistency of the estimator $\hat{F}^{(m, n)}$.

\begin{theorem} \label{thm:uniformconsistency}
Suppose that Assumption \ref{assump:uniformQhatconsistent} is satisfied. Then 
  \begin{align*}
  \|{F} - \hat{F}^{(m_n, n)}\|_{\mathcal{T}, \infty} \in \mathcal{O}_{\mathcal{P}_0}(g(n))
  \end{align*}
given that the equidistant grids 
$(\tau_k)_{k = 1}^{m_n}$ in $\mathcal{T}$ satisfy $\kappa_{m_n} \in o(g(n))$.
\end{theorem}

This shows that our estimator $\hat{F}^{(m, n)}$ can achieve uniform 
consistency over a set of distributions $\mathcal{P}_0 \subset \mathcal{P}$ 
given that the quantile regression procedure is uniformly consistent over 
$\mathcal{P}_0$. In Section \ref{sec:uniformlevelandpowerresults} we show how this strenghtened 
result can be used to establish asymptotic uniform level and power of our conditional 
independence test over $\mathcal{P}_0$. 

\subsection{A quantile regression model} \label{sec:quantileregmodels}

In this section we will provide an example of a flexible quantile regression model and 
estimation procedure where consistency results are available. Consider the model
  \begin{align} \label{eq:qrmodel}
  Q(\tau \mid z) = h(z)^T \beta_\tau 
  \end{align}
where $h : [0, 1]^d \to \R^p$ is a known and continuous transformation of $Z$, e.g., a polynomial or spline 
basis expansion to model non-linear effects. Inference in the model \eqref{eq:qrmodel} was 
analyzed by \cite{belloni2011} and \cite{belloni2019valid} in the high-dimensional setup $p \gg n$. 
In the following we describe a subset of their results that is relevant for our application.
Given an i.i.d. sample $(X_i, Z_i)_{i=1}^n$ and a fixed quantile regression level 
$\tau \in (0, 1)$, estimation of $\beta_\tau \in \R^p$ is 
carried out by penalized regression:
  \begin{align} \label{eq:betahat}
  \hat{\beta}_\tau \in \argmin_{\beta \in \R^p} \sum_{i = 1}^n L_\tau(X_i - h(Z_i)^T \beta) 
  + \lambda_\tau \|\beta\|_1
  \end{align}
where $L_\tau(u) = u(\tau - 1(u < 0))$ is the check function, $\|\cdot\|_1$ is the $1$-norm and 
$\lambda_\tau \geq 0$ is a tuning parameter that determines the degree of penalization. 
The tuning parameter $\lambda_\tau$ for a 
set $\mathcal{Q}$ of quantile regression levels can be chosen
in a data driven way as follows \citep[Section 2.3]{belloni2011}.  
Let $W_i = h(Z_i)$ denote the transformed predictors for $i=1, \dots, n$. Then we set 
  \begin{align} \label{eq:tuningparameter}
  \lambda_\tau = c \lambda \sqrt{\tau (1 - \tau)}
  \end{align}
where $c > 1$ is a constant with recommended value $c = 1.1$ and $\lambda$ is the $(1-n^{-1})$-quantile 
of the random variable 
  \begin{align*}
  \sup_{\tau \in \mathcal{T}} \frac{\| \Gamma^{-1}
  \frac{1}{n} \sum_{i=1}^n \left( \tau - 1(U_i \leq \tau) W_i \right) 
  \|_\infty}{\sqrt{\tau (1 - \tau)}}
  \end{align*}
where $U_1, \dots, U_n$ are i.i.d. $\mathcal{U}[0, 1]$. Here 
$\Gamma \in \R^{p \times p}$ is a diagonal matrix with $\Gamma_{kk} = \frac{1}{n} \sum_{i=1}^n (W_i)_k^2$.
The value of $\lambda$ is determined by simulation.


Sufficient regularity conditions under which the above estimation procedure can be proven to be
consistent are as follows.

\begin{myassump} \label{assump:highdimquantreg}

Denote by $f_{X \mid Z}$ the conditional density of $X$ given $Z$. Let $c > 0$ and $C > 0$ be constants.

\begin{itemize}
\item[(i)] There exists $s$ such that $\|\beta_\tau\|_0 \leq s$ 
for all $\tau \in \mathcal{Q} := [c, 1-c]$. 

\item[(ii)] $f_{X \mid Z}$ is continuously differentiable such that 
$f_{X \mid Z}(Q_{X \mid Z}(\tau \mid z) \mid z) \geq c$ for each $\tau \in \mathcal{Q}$ and 
$z \in [0, 1]^d$. Moreover, $\sup_{x \in [0, 1]} f_{X \mid Z}(x \mid z) \leq C$ and 
$\sup_{x \in [0, 1]} \partial_x f_{X \mid Z}(x \mid z) \leq C$. 

\item[(iii)] The transformed predictor $W = h(Z)$ satisfies $c \leq E((W^T \theta)^2) \leq C$ for 
all $\theta \in \R^p$ with $\|\theta\|_2 = 1$.
Moreover, $(E(\|W\|_\infty^{2q}))^{1/(2q)} \leq M_n$ for some $q > 2$ where $M_n$ satisfies
  \begin{align*}
  M_n^2 \leq \frac{\delta_n n^{1/2 - 1 / q}}{s \sqrt{\log(p \vee n)}}
  \end{align*}
and $\delta_n$ is some sequence tending to zero.
\end{itemize}
\end{myassump}

Assumption \ref{assump:highdimquantreg} (i) is a sparsity assumption, (ii) is a regularity condition on the conditional 
distribution, while (iii) is an assumption on the predictors. 
Examples of distributions for which Assumption \ref{assump:highdimquantreg} is 
satisfied are given in \cite{belloni2011} Section 2.5. This includes location 
models with Gaussian noise and location-scale models with bounded 
covariates, which in our setup with $Z \in [0, 1]^d$ means all location-scale 
models.

The following result 
\citep[Section 2.6]{belloni2011}
regarding the estimator $\hat{\beta}_{\tau}$ then holds.

\begin{theorem} \label{thm:belloni}

Assume that the tuning parameters $\{ \lambda_\tau \mid \tau \in \mathcal{Q} \}$ have been chosen according 
to \eqref{eq:tuningparameter}. Then
  \begin{align*}
  \sup_{\tau \in \mathcal{Q}} \|\beta_\tau - \hat{\beta}_\tau \|_2 \in 
  \mathcal{O}_P \left( 
  \sqrt{\frac{s \log(p \vee n)}{n}}
  \right) 
  \end{align*}
under Assumption \ref{assump:highdimquantreg}.
\end{theorem}

As a corollary of this consistency result we have the following.

\begin{corollary} \label{cor:Qhatconsistency}
Let $\hat{Q}(\tau \mid z) = h(z)^T \hat{\beta}_\tau$ be the predicted conditional quantile 
using the estimator $\hat{\beta}_{\tau}$. Then 
  \begin{align*}
  \sup_{z \in [0, 1]^d} \sup_{\tau \in \mathcal{Q}} |Q(\tau \mid z) - \hat{Q}^{(n)}(\tau \mid z)|
  \in \mathcal{O}_P \left( 
  \sqrt{\frac{s \log(p \vee n)}{n}}
  \right) 
  \end{align*}
under Assumption \ref{assump:highdimquantreg}.
\end{corollary}

This shows that Assumption \ref{assump:Qhatconsistency} is satisfied 
under the model \eqref{eq:qrmodel} 
whenever Assumption \ref{assump:highdimquantreg} is satisfied 
with $\mathcal{T} \subset \mathcal{Q}$ and 
$\sqrt{s \log(p \vee n) / n} \to 0$, which is the key underlying assumption of 
Theorem \ref{thm:pointwiseconsistency}. Note also that 
Assumption \ref{assump:Qhatconsistency} (ii) is contained in 
Assumption \ref{assump:highdimquantreg} (ii). 
Theorem \ref{thm:belloni} and Corollary \ref{cor:Qhatconsistency} can 
be extended to hold uniformly over $\mathcal{P}_0 \subset \mathcal{P}$ 
by assuming that the conditions of Assumption \ref{assump:highdimquantreg} 
hold uniformly over $\mathcal{P}_0$. This then gives the 
statement of Assumption \ref{assump:uniformQhatconsistent} that is 
required for Theorem \ref{thm:uniformconsistency}. 

\section{Testing Conditional Independence} \label{sec:test-cond-ind}

In this section we describe the conditional independence testing 
framework in terms of the so-called partial copula. 
As above we let $(X, Y, Z) \sim P \in \mathcal{P}$ such that $X, Y \in [0, 1]$ and $Z \in [0, 1]^d$ where $\mathcal{P}$ are the distributions 
that are absolutely continuous with respect to Lebesgue measure 
on $[0, 1]^{2 + d}$. 
Also let $f$ denote a generic density function. We then say that 
$X$ is conditionally independent of $Y$ given $Z$ if 
  \begin{align*}
  f(x, y \mid z) = f(x \mid z) f(y \mid z)
  \end{align*}
for almost all $x, y \in [0, 1]$ and $z \in [0, 1]^d$. See e.g. \cite{dawid1979conditional}.
In this case we write that
$X \indep_P Y \mid Z$, where we usually omit the dependence on $P$ when there is no ambiguity.
By $\mathcal{H} \subset \mathcal{P}$ we denote the subset of distributions 
for which conditional independence is satisfied, and we let 
$\mathcal{Q} := \mathcal{P} \setminus \mathcal{H}$ be 
the alternative of conditional dependence.

\subsection{The partial copula} \label{sec:pit}

We can regard the conditional distribution function as a mapping 
$(t, z) \mapsto F(t \mid z)$ for $t \in [0, 1]$ and $z \in [0, 1]^d$. Assuming
that this mapping is measurable, we define a new pair of random variables 
$U_1$ and $U_2$ by the transformations
  \begin{align*}
  U_1 := F_{X \mid Z}(X \mid Z) \quad \text{and} \quad 
  U_2 := F_{Y \mid Z}(Y \mid Z).
  \end{align*}
This transformation is usually called the probability integral transformation or 
Rosenblatt transformation due to \cite{rosenblatt1952remarks}, where the transformation 
was initially introduced and the following key result was shown.

\begin{proposition} \label{prop:Uuniform}
It holds that $U_\ell \sim \mathcal{U}[0, 1]$ and $U_\ell \indep Z$ for $\ell=1, 2$.
\end{proposition}

Hence the transformation can be understood as a normalization, where marginal dependencies
of $X$ on $Z$ and $Y$ on $Z$ are filtered away. 
The joint distribution of $U_1$ and $U_2$ has been termed 
the partial copula of $X$ and $Y$ given $Z$ in the copula literature 
\citep{bergsma_nonparametric_2011, spanhel_partial_2016}.
Independence in the partial copula relates 
to conditional independence in the following way.

\begin{proposition} \label{prop:XYZtoUU}
If $X \indep Y \mid Z$ then $U_1 \indep U_2$. 
\end{proposition}

Therefore the question about conditional independence can be transformed into a 
question about independence. Note, however, that $U_1 \indep U_2$ does not in general 
imply that $X \indep Y \mid Z$. See Property 7 in 
\cite{spanhel_partial_2016} for a counterexample

The variables $U_\ell$ were termed nonparametric residuals by \cite{patra2016nonparametric}
due to the independence property $U_\ell \indep Z$ 
which is analogues to the property of 
conventional residuals in additive Gaussian noise models. Note that the entire 
conditional distribution function is required in order to compute the 
nonparametric residual, while conventional residuals in additive noise models 
are computed using only the conditional expectation. 
In return, Proposition \ref{prop:Uuniform} provides the distribution of 
the nonparametric residuals without asumming any functional or distributional 
relationship between $X$ ($Y$ resp.) and $Z$, whereas the distribution 
of conventional residuals is not known without further assumptions. 
Moreover, the nonparametric residuals $U_1$ and $U_2$ 
are independent under conditional independence, while
conventional residuals are only uncorrelated 
unless we make a Gaussian assumption, say.

\subsection{Generic testing procedure} \label{sec:testingprocedure}

Suppose $(X_i, Y_i, Z_i)_{i=1}^n$ is a sample from $P \in \mathcal{P}_0$ where 
$\mathcal{P}_0$ is some subset of $\mathcal{P}$. 
Also let $\mathcal{H}_0 := \mathcal{P}_0 \cap \mathcal{H}$ 
and $\mathcal{Q}_0 := \mathcal{P}_0 \cap \mathcal{Q}$ be the distributions in 
$\mathcal{P}_0$ satisfying conditional independence and 
conditional dependence, respectively. Denote by
  \begin{align*}
  U_{1, i} := F_{X \mid Z}(X_i \mid Z_i) \quad
  \text{and} \quad
  U_{2, i} := F_{Y \mid Z}(Y_i \mid Z_i)
  \end{align*}
the nonparametric residuals for $i=1, \dots, n$.
Let $\psi_n : [0, 1]^{2n}\to \{0, 1\}$ denote a test for 
independence in a bivariate continuous distribution. The observed value of the test is 
  \begin{align*}
  \Psi_n := \psi_n((U_{1, i}, U_{2, i})_{i=1}^n)
  \end{align*}
with $\Psi_n = 0$ indicating acceptance and $\Psi_n = 1$ rejection of the hypothesis. 
By Proposition \ref{prop:XYZtoUU} we then reject the hypothesis of 
conditional independence, $X \indep Y \mid Z$,  if $\Psi_n = 1$. 
However, in order to implement
the test in practice, we will need to replace the conditional distribution functions 
$F_{X \mid Z}$ and $F_{Y \mid Z}$ by estimates.

Given some generic estimators of the conditional distribution functions we can formulate 
a generic version of the partial copula conditional independence test as follows. 

\pagebreak

\begin{definition} \label{def:naivetest}
Let $(X_i, Y_i, Z_i)_{i=1}^n$ be an i.i.d. sample from $P \in \mathcal{P}_0$. Also let 
$\psi_n$ be a test for independence in a bivariate continuous distribution.
  \begin{itemize}
  \item[(i)] Form the estimates $\hat{F}^{(n)}_{X \mid Z}$ and $\hat{F}^{(n)}_{Y \mid Z}$ based on 
  $(X_i, Y_i, Z_i)_{i=1}^n$.
  
  \item[(ii)] Compute the estimated nonparametric residuals
    \begin{align*}
    \hat{U}^{(n)}_{1, i} := \hat{F}^{(n)}_{X \mid Z}(X_i \mid Z_i) \quad \text{and} 
    \quad \hat{U}^{(n)}_{2, i} := \hat{F}^{(n)}_{Y \mid Z}(Y_i \mid Z_i)
    \end{align*}
  for $i = 1, \dots, n$.
  
  \item[(iii)] Let 
  $\hat{\Psi}_n : = \psi_n \left( (\hat{U}_{1, i}^{(n)}, \hat{U}_{2, i}^{(n)})_{i=1}^n \right)$ 
  and reject the hypothesis $X \indep Y \mid Z$ of conditional independence
  if $\hat{\Psi}_n = 1$. 
  \end{itemize}
\end{definition}

This generic version of the conditional independence 
test is analogous to the approach of 
\cite{bergsma_nonparametric_2011}, but here we emphasize the 
modularity of the testing procedure. Firstly, one can 
use any method for estimating conditional distribution functions. Secondly, 
any test for independence in a bivariate continuous distribution can be utilized.

We note that under the assumptions of Theorem \ref{thm:pointwiseconsistency} 
it holds that   
	\begin{align*}
	|(\hat{U}^{(n)}_{1,i}, \hat{U}^{(n)}_{2,i}) - 
	(U_{1,i}, U_{2,i})|_{\mathcal{T}, 1} \stackrel{P}{\rightarrow} 0
	\end{align*}
where $|u - v|_{\mathcal{T}, 1} = |u_1 - v_1|1(u_1, v_1 \in \mathcal{T}) + 
|u_2 - v_2|1(u_2, v_2 \in \mathcal{T})$. That is, each estimated pair of nonparametric residuals has the partial copula as asymptotic distribution -- 
except perhaps on the fringe part of the unit square outside of 
$\mathcal{T} \times \mathcal{T}$. This is a priori only a marginal 
result for each $i$, but it suggests that tests based on the estimated 
residuals behave as if they were i.i.d. observations from the partial copula. 

Once we have chosen the test for independence, $\psi_n$, we can establish rigorous results on the properties of the test over the space of hypotheses $\mathcal{H}_0$ and alternatives $\mathcal{Q}_0$, but how exactly to transfer the consistency of the estimated residuals to results on level and power depends on the specific test statistic. We will in the following sections demonstrate this transfer for one particular class of test statistics.

\subsection{Generalized measure of correlation} \label{sec:gcm}

We will now introduce a generalized measure of correlation that will form the basis 
for an independence test between the nonparametric residuals 
$U_1$ and $U_2$.

\begin{definition}
The generalized correlation, $\rho$, between 
$U_1$ and $U_2$ is defined in term of a multivariate 
function $\phi = (\phi_1, \dots, \phi_q) : [0, 1] \to \R^q$ as
	\begin{align}
	\rho = E_P(\phi(U_1) \phi(U_2)^T)
	\end{align}
such that $\rho$ is a $q \times q$ matrix with entries 
$\rho_{k \ell} = E_P(\phi_k(U_1) \phi_\ell(U_2))$ 
for $k, \ell = 1, \dots, q$.
\end{definition}

We will assume that the function $\phi = (\phi_1, \dots, \phi_q)$ 
defining the generalized correlation 
satisfies the following assumptions for the remainder of the paper.
\pagebreak

	\begin{myassump} \label{assump:phimap}
	\
		\begin{itemize}
		\item[(i)] The support $\mathcal{T}_k$ of each coordinate function 
		$\phi_k$ is a compact subset of $(0, 1)$. 
		
		\item[(ii)] Each coordinate function $\phi_k : [0, 1] \to \R$ is 
		Lipschitz continuous.
		
		\item[(iii)] $\int_0^1 \phi_k(u) \mathrm{d} u = 0$ and 
		$\int_0^1 \phi_k(u)^2 \mathrm{d} u = 1$ for each $k=1, \dots, q$.
		
    \item[(iv)] The coordinate functions $\phi_1, \dots, \phi_q$ are linearly
    independent.
		\end{itemize}
	\end{myassump}
Let us provide some intuition about the interpretation of the 
generalized correlation $\rho$ and explain the role of the assumptions 
on $\phi$ in Assumption \ref{assump:phimap}.

Each entry $\rho_{k \ell}$ can be interpreted as an expected conditional 
correlation, and it can be understood in terms of the partial 
and conditional copula \citep{patton2006}. Let 
$C(u_1, u_2 \mid z) = F(U_1 \leq u_1, U_2 \leq u_2 \mid Z=z)$ denote the conditional copula of 
$X$ and $Y$ given $Z=z$. Then the partial copula is the 
expected conditional copula, i.e., $C_p(u_1, u_2) = E_P(C(u_1, u_2 \mid Z))$ \citep{spanhel_partial_2016}.
The conditional generalized correlation, $\rho_{k \ell}(z)$,
between $X$ and $Y$ given $Z=z$ can be expressed in terms of the 
conditional copula by
  \begin{align*}
  \rho_{k \ell}(z) := E_P(\phi_k(U_1) \phi_\ell(U_2) \mid Z = z) = \int \phi_k(u_1) \phi_\ell(u_2) C(\mathrm{d}u_1 , \mathrm{d}u_2 \mid z).
  \end{align*}
By the tower property of conditional expectations, $\rho_{k \ell}$ can be represented as an 
expected generalized correlation  
  \begin{align*}
  \rho_{k \ell} = 
  E_P(\rho_{k \ell}(Z)) = \int \phi_k(u_1) \phi_\ell(u_2) 
  C_p(\mathrm{d}u_1, \mathrm{d}u_2).
  \end{align*}
Hence $\rho_{k \ell}$ measures the expected conditional generalized correlation of 
$X$ and $Y$ given $Z$ w.r.t. the distribution of $Z$. Importantly, Assumption 
\ref{assump:phimap} (iii) implies that
	\begin{align*}
	\rho = E_P(\phi(U_1) \phi(U_2)^T) = E_P(\phi(U_1)) E(\phi(U_2)^T) = 0
	\end{align*}
whenever $X \indep Y \mid Z$ due to Proposition \ref{prop:XYZtoUU}.

\begin{figure}[t]
\centering
\includegraphics[scale=0.7]{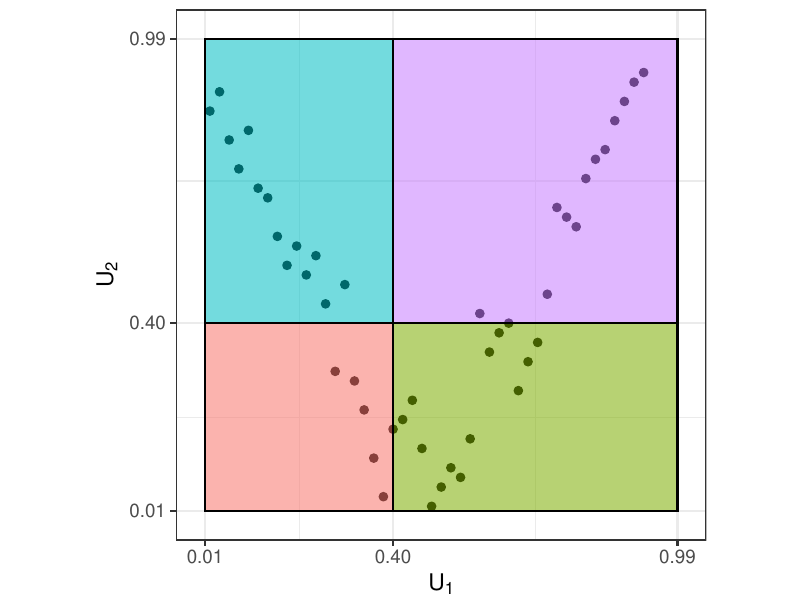}
\caption{A sample from a copula $(U_1, U_2)$ with clear dependence, but where the overall sample correlation is close to zero. The dependence is captured 
by considering sample correlation of observations in different regions of 
the outcome space.}
\label{fig:rhoplot}
\end{figure}

The purpose of Assumption \ref{assump:phimap} (i) is twofold. Firstly, letting 
the supports $\mathcal{T}_k$ and $\mathcal{T}_\ell$ of $\phi_k$ and $\phi_\ell$ be subsets of 
$(0, 1)$ implies that $\rho_{k \ell}$ focuses on dependence in the compact region 
$\mathcal{T}_k \times \mathcal{T}_\ell \subset (0, 1)^2$ of the outcome space $[0, 1]^2$ of 
$(U_1, U_2)$. For $q \geq 2$ the generalized correlation 
$\rho$ thus summarizes dependencies in different regions of the outcome space. 
See Figure \ref{fig:rhoplot} for an illustration of this idea. Secondly, 
the supports $(\mathcal{T}_k)_{k=1}^q$ will play the role as subsets of the possible quantile levels 
$\mathcal{T}  = [\tau_{\min}, \tau_{\max}]$, when we choose the conditional distribution 
function estimators based on quantile regression from Section \ref{sec:quantregbasedestimator}. 
This connection will be made clear in Section \ref{sec:correlationtest}.

The functional form of $\phi_k$ and $\phi_\ell$ determines the kind of 
dependence measured by $\rho_{k \ell}$. Ignoring Assumption \ref{assump:phimap} (i), 
consider letting $\phi_k(u) = \phi_\ell(u) = \sqrt{12}(u - 1 / 2)$ for $u \in [0, 1]$. Then 
$\rho_{k \ell}$ measures the expected conditional Spearman correlation between $X$ and $Y$ given 
$Z$ with respect to the distribution of $Z$. In Section \ref{sec:trimming} we describe a choice 
of functions $\phi_k$ that leads to a trimmed version of expected conditional Spearman correlation which 
satisfies Assumption \ref{assump:phimap} (i). As we shall see in Section \ref{sec:correlationtest}, 
Assumption \ref{assump:phimap} (ii), i.e., that the coordinate 
functions $\phi_k$ are Lipschitz continuous, is 
crucial for deriving asymptotic properties for the empirical version of the generalized correlation $\rho$.
Lastly, we assume that $\phi_1, \dots, \phi_q$ are linearly independent 
in Assumption \ref{assump:phimap} (iv) to 
avoid degeneracy of its empirical version, which we introduce in Section 
\ref{sec:correlationtest}.

\subsection{Test based on generalized correlation} \label{sec:correlationtest}

In this section we will analyze in depth the conditional independence 
test resulting from basing the test $\psi_n$ in Definition 
\ref{def:naivetest} on the generalized correlation $\rho$. 
We will formulate the results in terms of a 
generic method for estimating conditional distribution functions in order to emphasize 
the generality of the method and 
illustrate the abstract assumptions needed for the 
test to be sound. Along the way we will explain when the assumptions are satisfied for 
the quantile regression based estimator $\hat{F}^{(m, n)}$ that we developed in Section 
\ref{sec:conddistest}.

With $\rho$ the generalized correlation between $U_1$ and $U_2$ defined in terms of a function 
$\phi$ satisfying Assumption \ref{assump:phimap}
we let $\rho_n : [0, 1]^{2n} \to \R^{q \times q}$ be its corresponding empirical version:
	\begin{align} \label{eq:rhoteststatistic}
	\rho_n((u_i, v_i)_{i=1}^n) := \frac{1}{n} \sum_{i=1}^n \phi(u_i) \phi(v_i)^T.
	\end{align}
Soundness of a test based on $\rho_n$ depends on consistency of the estimators 
$\hat{F}_{X \mid Z}^{(n)}$ and $\hat{F}_{Y \mid Z}^{(n)}$.
Recall that we by $\mathcal{T}_1, \dots, \mathcal{T}_q$ 
denote the supports of $\phi_1, \dots, \phi_q$.
Let $\tau_{\min} := \inf (\mathcal{T}_1 \cup 
\cdots \cup \mathcal{T}_q) > 0$ 
and $\tau_{\max} := 
\sup (\mathcal{T}_1 \cup \cdots \cup \mathcal{T}_q) < 1$, and then define
$\mathcal{T} := [\tau_{\min}, \tau_{\max}]$. 
As in Section \ref{sec:quantregbasedestimator} we let the norm
$\|\cdot \|_{\mathcal{T}, \infty}$ be given by
	\begin{align*}
	\| f(t, z) \|_{\mathcal{T}, \infty} = \sup_{z \in [0, 1]^d} 
	\sup_{t \in Q_{X \mid Z}(\mathcal{T} \mid z)} |f(t, z)|
	\end{align*}
when $X$ given $Z$ is the conditional distribution of interest. Similarly define 
$\|\cdot \|_{\mathcal{T}, \infty}'$ by 
	\begin{align*}
	\| f(t, z) \|_{\mathcal{T}, \infty}' = \sup_{z \in [0, 1]^d} 
	\sup_{t \in Q_{Y \mid Z}(\mathcal{T} \mid z)} |f(t, z)|.
	\end{align*}
Then we have the following assumption on our estimators.

\begin{myassump} \label{assump:Fhatconsistency}
For each distribution $P \in \mathcal{P}_0$ there exist deterministic rate functions 
$g_P$ and $h_P$ tending to zero as $n \to \infty$ and functions 
$\xi, \xi' : [0, 1] \times [0, 1]^d \to \R$ such that
	\begin{itemize}
	\item[(i)] 
	$\| F_{X \mid Z} - \hat{F}_{X \mid Z}^{(n)}\|_{\mathcal{T}, \infty} 
	\in \mathcal{O}_P(g_P(n))$ and 
	$\| F_{Y \mid Z} - \hat{F}_{Y \mid Z}^{(n)}\|_{\mathcal{T}, \infty}'
	\in \mathcal{O}_P(h_P(n))$.
	
	\item[(ii)] $||\xi - \hat{F}^{(n)}_{X \mid Z}||_{\mathcal{T}^c, 
	\infty} \in \mathcal{O}_P(g_P(n))$ and 
	$||\xi' - \hat{F}^{(n)}_{Y \mid Z}||_{\mathcal{T}^c, \infty} '
	\in \mathcal{O}_P(h_P(n))$.
	\end{itemize}
\end{myassump}

Assumption \ref{assump:Fhatconsistency} (i)
states that our estimators 
$\hat{F}^{(n)}_{X \mid Z}$ and $\hat{F}^{(n)}_{Y \mid Z}$
are consistent with rates $g_P$ and $h_P$ over 
the conditional $\mathcal{T}$-quantiles in their respective 
conditional distributions. 
This is the result of Theorem \ref{thm:pointwiseconsistency} 
regarding our quantile regression based estimator 
$\hat{F}^{(m, n)}$ when $\mathcal{T}$ as above is taken as the set of
potential quantile regression levels.

Assumption \ref{assump:Fhatconsistency} (ii) is an 
assumption on the behavior of our estimator in the 
tails of the conditional distribution, i.e., over the 
conditional $\mathcal{T}^c$-quantiles. Here we do not assume 
consistency, but we do assume that the limit for $n \to \infty$ exists, 
and that our estimators are convergent to their limits with 
rates $g_P$ and $h_P$ respectively. This assumption 
is satisfied by our quantile regression based estimator 
$\hat{F}^{(m, n)}$ whenever it satisfies 
Assumption \ref{assump:Fhatconsistency} (i).

With this assumption we first establish the 
asymptotic distribution of the test statistic
	\begin{align} \label{eq:rhohat}
	\hat{\rho}_n := \rho_n \left(
	(\hat{U}_{1, i}, \hat{U}_{2, i})_{i=1}^n
	\right) = 
	\frac{1}{n} \sum_{i=1}^n \phi(\hat{U}_{1, i}) 
	\phi(\hat{U}_{2, i})^T
	\end{align}
under the hypothesis of conditional independence. Below we use 
$\Rightarrow_P$ to denote convergence in distribution with respect to $P$.

\begin{theorem} \label{thm:rhoasympnormal}
Suppose that Assumption \ref{assump:Fhatconsistency} is satisfied with rate functions 
$g_P$ and $h_P$ such that 
$\sqrt{n} g_P(n) h_P(n) \to 0$ as $n \to \infty$ for each $P \in \mathcal{P}_0$. 
Then the statistic $\hat{\rho}_n$ given by 
\eqref{eq:rhohat} satisfies 
  \begin{align*}
  \sqrt{n} \hat{\rho}_n \Rightarrow_P \mathcal{N}(0, \Sigma \otimes \Sigma)
\end{align*}   
for each fixed $P \in \mathcal{H}_0$. The asymptotic covariance matrix
is given by
	\begin{align*}
	\Sigma_{k, s} = \int_0^1 \phi_k(u) \phi_s(u) \mathrm{d} u 
	\end{align*}
for $k, s = 1, \dots, q$ and does not depend on $P$.
\end{theorem} 

If the rate functions are $g_P(n) = n^{-a}$ and $h_P(n) = n^{-b}$, 
then we require that $a + b > 1/2$. Thus
convergence slightly faster than rate $n^{-1/4}$ for both 
estimators is sufficient, but there can 
also be a tradeoff between the rates. 
Interestingly, Theorem \ref{thm:rhoasympnormal} does 
not require sample splitting for the estimation of 
the conditional distribution function and 
computation of the test statistic. This is due to 
the fact that we are only interested in the asymptotic distribution 
under conditional independence. A similar 
phenomenon was found by \cite{shah2018hardness}, 
when they proved asymptotic normality of 
their Generalised Covariance Measure under conditional independence.

According to Corollary \ref{cor:Qhatconsistency},
the requirement $\sqrt{n} g_P(n) h_P(n) \to 0$ is satisfied 
for our quantile regression based estimator $\hat{F}^{(m, n)}$, if 
the quantile regression model \eqref{eq:qrmodel} of Section 
\ref{sec:quantileregmodels} is valid for both 
$X$ given $Z$ and $Y$ given $Z$ for some continuous transformations 
$h_1 : [0, 1]^d \to \R^{p_1}$ and $h_2 : [0, 1]^d \to \R^{p_2}$ 
and Assumption 
\ref{assump:highdimquantreg} is satisfied with $s_1, s_2, p_1, p_2, n \to \infty$ 
such that
	\begin{align*}
	\sqrt{n} \cdot \sqrt{\frac{s_1 \log (p_1 \vee n)}{n}} 
	\cdot \sqrt{\frac{s_2 \log (p_2 \vee n)}{n}} 
	= 
	\sqrt{\frac{s_1 s_2 \log (p_1 \vee n) \log(p_2 \vee n)}{n}} 
	\to 0
	\end{align*}
where $s_1 = \sup_{\tau \in \mathcal{T}} \|\beta_{1, \tau} \|_0$ and 
$s_2 = \sup_{\tau \in \mathcal{T}} \|\beta_{2, \tau} \|_0$ are the 
sparsities of the model parameters.

With the test statistic
	\begin{align} \label{eq:chistatistic}
	T_n := \| \Sigma^{-1 / 2} \hat{\rho}_n \Sigma^{-1 / 2} \|_F^2,
	\end{align}
where $\| \cdot \|_F$ denotes the Frobenius norm, we have the following corollary of Theorem \ref{thm:rhoasympnormal}.

\begin{corollary} \label{cor:chisqteststat}
Let the condition of Theorem \ref{thm:rhoasympnormal} be satisfied 
and let $T_n$ be given by \eqref{eq:chistatistic}. Then it holds that
	\begin{align*}
	n T_n \Rightarrow_P \chi^2_{q^2}
	\end{align*}
for each fixed $P \in \mathcal{H}_0$. 
\end{corollary}

In view of Theorem \ref{thm:rhoasympnormal} and Corollary \ref{cor:chisqteststat} we define the following
conditional independence test based on the generalized correlation.

\begin{definition} \label{def:correlationtest}
Let $\alpha \in (0, 1)$ be a desired significance level and 
$T_n$ the test statistic \eqref{eq:chistatistic}. Then we let 
$\hat{\Psi}_n$ be the test given by
	\begin{align*}
	\hat{\Psi}_n = 1(T_n > 
	n^{-1} z_{1 - \alpha})
	\end{align*}
where $z_{1 - \alpha}$ is the $(1-\alpha)$-quantile of a 
$\chi^2_{q^2}$-distribution.
\end{definition}

Control of the asymptotic pointwise level is then an easy corollary of Corollary \ref{cor:chisqteststat}.

\begin{corollary} \label{cor:level}
Suppose that Assumption \ref{assump:Fhatconsistency} is satisfied with rate functions 
$g_P$ and $h_P$ such that 
$\sqrt{n} g_P(n) h_P(n) \to 0$ as $n \to \infty$ for each $P \in \mathcal{P}_0$.  Then the 
test $\hat{\Psi}_n$ given by 
Definition \ref{def:correlationtest} has asymptotic pointwise level over $\mathcal{H}_0$, 
i.e.,
  \begin{align*}
  \limsup_{n \to \infty} E_P(\hat{\Psi}_n) = \alpha
  \end{align*}
for each fixed $P \in \mathcal{H}_0$. 
\end{corollary}

This shows that the test achieves correct 
level given consistency of the estimators 
$\hat{F}^{(n)}_{X \mid Z}$ and $\hat{F}^{(n)}_{Y \mid Z}$ with suitably fast rates. To
obtain results on power of the test $\hat{\Psi}_n$ 
we only need to understand how $\hat{\rho}_n$ converges in probability
and not its entire asymptotic distribution.

\begin{theorem} \label{thm:rhohat_to_rho}

The test statistic $\hat{\rho}_n$ given by 
\eqref{eq:rhohat} satisfies
  \begin{align*}
  \hat{\rho}_n \stackrel{P}{\longrightarrow} \rho 
  \end{align*}
for each fixed $P \in \mathcal{P}_0$ under
Assumption \ref{assump:Fhatconsistency}.
\end{theorem}

One may note that the theorem does not require that the rate functions $g_P$ 
and $h_P$ converge to zero at a certain rate.
Let $\mathcal{A}_0 \subseteq \mathcal{Q}_0$ be the subset of 
alternatives for which $\rho_{k \ell} \neq 0$ for at least one 
combination of $k, \ell = 1, \dots, q$. 
Then we have the following corollary 
of Theorem \ref{thm:rhohat_to_rho}, which exploits that $n T_n$
diverges in probability whenever $P \in \mathcal{A}_0$. 

\begin{corollary} \label{cor:power}
For each level $\alpha \in (0, 1)$ the test $\hat{\Psi}_n$ given by 
Definition \ref{def:correlationtest} has asymptotic pointwise power against $\mathcal{A}_0$, i.e.,
  \begin{align*}
   \liminf_{n \to \infty} E_P(\hat{\Psi}_n) = 1
  \end{align*}
for each fixed $P \in \mathcal{A}_0$ under Assumption \ref{assump:Fhatconsistency}.
\end{corollary}

Let us discuss the alternatives the test has power against. 
Firstly, note that we always have the implications
	\begin{align*}
	X \indep Y \mid Z \quad \Rightarrow \quad 
	U_1 \indep U_2 \quad \Rightarrow \quad \rho = 0
	\end{align*}
However, none of the reverse implications 
are in general true. We do, however, have the 
following result stating a sufficient condition for the reverse implication 
of the first statement. 

\begin{proposition} \label{prop:UUtoXYZ}
Assume that $(U_1, U_2) \indep Z$. 
Then $X \indep Y \mid Z$ if and only if $U_1 \indep U_2$.
\end{proposition}

This means that if $Z$ only affects the marginal 
distributions of $X$ and $Y$, then independence in the partial copula 
implies conditional independence. This is known as the simplifying 
assumption in the copula literature \citep{gijbels_estimation_2015, 
spanhel2015partial}. Naturally, $U_1 \notindep U_2$ always implies 
$X \notindep Y \mid Z$, so the simplifying assumption is not a necessary 
condition for our test to have power, but it does give some intuition about a 
subset of distributions for which the partial copula completely characterizes 
conditional independence. However, an unavoidable 
limitation of the method is that it can never have power against 
alternatives for which $U_1 \indep U_2$ but $X \notindep Y \mid Z$. 

Turning to the second implication, Corollary \ref{cor:power} tells us that 
we have power against alternatives for which 
$\rho_{k \ell} \neq 0$ for some $k, \ell = 1, \dots, q$.
However, not all types of dependencies can be detected in this fashion, and 
it is possible that $\rho = 0$, while $U_1 \notindep U_2$. A test based on 
$\rho$ will not have power against such an alternative. For an abstract 
interpretation of the generalized correlation $\rho$ we refer to Section 
\ref{sec:gcm}. In Section \ref{sec:trimming} we introduce a concrete 
generalized correlation and elaborate on its interpretation.

Finally, basing the test on values of $T_n$ is natural
since the asymptotic behaviour is readily available through 
Theorem \ref{thm:rhoasympnormal}, but other transformations 
of $\hat{\rho}_n$ could be considered such as taking the coordinatewise 
absolute maximum $\max_{k,l} | (\Sigma^{-1 / 2} \hat{\rho}_n \Sigma^{-1 / 2})_{k,l} | = 
\| \Sigma^{-1 / 2} \hat{\rho}_n \Sigma^{-1 / 2} \|_{\infty}$.

\subsection{Uniform level and power results} \label{sec:uniformlevelandpowerresults}

The level and power results of Section \ref{sec:correlationtest}
are pointwise over 
the space of hypotheses and alternatives, i.e., they state level and power of the test 
when fixing a distribution $P$. 
In this section we describe how these results 
can be extended to hold uniformly by 
strengthening the statements in Assumption \ref{assump:Fhatconsistency} to hold uniformly.

\begin{myassump} \label{assump:uniformFhatconsistency}

For $\mathcal{P}_0 \subset \mathcal{P}$ 
there exist deterministic rate functions 
$g$ and $h$ tending to zero as $n \to \infty$ and functions 
$\xi, \xi' : [0, 1] \times [0, 1]^d \to \R$ such that
	\begin{itemize}
	\item[(i)] 
	$\| F_{X \mid Z} - \hat{F}_{X \mid Z}^{(n)}\|_{\mathcal{T}, \infty} 
	\in \mathcal{O}_{\mathcal{P}_0}(g(n))$ and 
	$\| F_{Y \mid Z} - \hat{F}_{Y \mid Z}^{(n)}\|_{\mathcal{T}, \infty}'
	\in \mathcal{O}_{\mathcal{P}_0}(h(n))$.
	
	\item[(ii)] $||\xi - \hat{F}^{(n)}_{X \mid Z}||_{\mathcal{T}^c, 
	\infty} \in \mathcal{O}_{\mathcal{P}_0}(g(n))$ and 
	$||\xi' - \hat{F}^{(n)}_{Y \mid Z}||_{\mathcal{T}^c, \infty} '
	\in \mathcal{O}_{\mathcal{P}_0}(h(n))$.
	\end{itemize}
\end{myassump}

As before we note that Assumption \ref{assump:uniformFhatconsistency} (i) 
is the result of Theorem \ref{thm:uniformconsistency} regarding 
our quantile regression based estimator $\hat{F}^{(m, n)}$. Moreover, 
Assumption \ref{assump:uniformFhatconsistency} (ii) is valid for 
$\hat{F}^{(m, n)}$ whenever it satisfies Assumption \ref{assump:uniformFhatconsistency} (i). 

We will now describe the extensions of Theorem \ref{thm:rhoasympnormal} 
and Theorem \ref{thm:rhohat_to_rho} that can be 
obtained under Assumption \ref{assump:uniformFhatconsistency}. Below we 
write $\Rightarrow_\mathcal{M}$ to denote uniform convergence in distribution over a 
set of distributions $\mathcal{M}$, and we use $\to_\mathcal{M}$ to denote uniform 
convergence in probability over $\mathcal{M}$.
We refer to Appendix \ref{appendix:modesofconvergence}
for the formal definitions.

\begin{theorem} \label{thm:uniformresults}

Let $\hat{\rho}_n$ be the statistic given by \eqref{eq:rhohat}. Then we have:

\begin{itemize}
\item[(i)] Under Assumption \ref{assump:uniformFhatconsistency} 
with rate functions satisfying $\sqrt{n} g(n) h(n) \to 0$ 
it holds that  
  \begin{align*}
  \sqrt{n} \hat{\rho}_n \Rightarrow_{\mathcal{H}_0} \mathcal{N}(0, \Sigma \otimes \Sigma)
  \end{align*}
where $\Sigma$ is as in Theorem 
\ref{thm:rhoasympnormal}.

\item[(ii)] Under Assumption \ref{assump:uniformFhatconsistency} it holds that $\hat{\rho}_n \to_{\mathcal{P}_0} \rho$. 
\end{itemize}
\end{theorem}

As a straightforward corollary of Theorem \ref{thm:uniformresults} (i) 
we get the following uniform level result. 

\begin{corollary} \label{cor:uniformlevel}
The test $\hat{\Psi}_n$ given by 
Definition \ref{def:correlationtest} has asymptotic uniform level over $\mathcal{H}_0$, i.e.,
  \begin{align*}
  \limsup_{n \to \infty} \sup_{P \in \mathcal{H}_0} E_P(\hat{\Psi}_n) = \alpha,
  \end{align*}
given that Assumption \ref{assump:uniformFhatconsistency} is satisfied with 
$\sqrt{n} g(n) h(n) \to 0$ as $n \to \infty$. 
\end{corollary}

The pointwise power result of Theorem \ref{cor:power} does not extend directly to 
a uniform version in the same way as the level result. 
For $\lambda > 0$ we let $\mathcal{A}_\lambda \subset \mathcal{Q}_0$ 
be the set of alternatives for which $|(\rho_P)_{k \ell}| > \lambda$ for at 
least one combination of $k, \ell = 1, \dots, q$, where we emphasize 
that $\rho_P$ depends on the distribution $P$.
We then have the following uniform power result as a corollary of 
Theorem \ref{thm:uniformresults} (ii).

\begin{corollary} \label{cor:uniformpower}
For all fixed levels $\alpha \in (0, 1)$ the test $\hat{\Psi}_n$ given by 
Definition \ref{def:correlationtest} has asymptotic uniform power against 
$\mathcal{A}_{\lambda}$ for each $\lambda > 0$, i.e.,
  \begin{align*}
   \liminf_{n \to \infty} \inf_{P \in \mathcal{A}_\lambda} E_P(\hat{\Psi}_n) = 1,
  \end{align*}
under Assumption \ref{assump:uniformFhatconsistency}.
\end{corollary}

The reason we need to restrict to the class of alternatives $\mathcal{A}_\lambda$ for 
a fixed $\lambda > 0$ is the 
following. If the infimum is taken over $\mathcal{A}_0$, then there could exist a 
sequence $(P_m)_{m=1}^\infty \subset \mathcal{A}_0$ of distributions such that $(\rho_{P_m})_{k \ell} \neq 0$ 
for each $m \geq 1$ but $(\rho_{P_m})_{k \ell} \to 0$ as $m \to \infty$. As a consequence 
$n T_n$ will not necessarily diverge in probability, which is 
crucial to the proof of the corollary. However, when restricting to $\mathcal{A}_\lambda$ 
we are ensured that $\inf_{P \in \mathcal{A}_\lambda} 
|(\rho_P)_{k \ell}| \geq \lambda > 0$ for at least one combination of 
$k, \ell = 1, \dots, q$. 

We note that these uniform level and power results are not in 
contradiction with the impossibility result of 
\cite{shah2018hardness} mentioned in Section \ref{sec:introduction} 
because our results apply to a restricted set of distributions, 
$\mathcal{P}_0$, where the conditional distribution functions 
are estimable with sufficiently fast rate.

\subsection{Trimmed Spearman correlation} \label{sec:trimming}

We will now define a specific family of functions $\phi$ defining 
the generalized correlation that can be shown to satisfy 
Assumption \ref{assump:phimap}, which results in trimmed versions 
of the expected conditional Spearman correlation. As mentioned 
in Section \ref{sec:gcm}, ignoring Assumption \ref{assump:phimap} (i), 
we could consider
	\begin{align} \label{eq:spearman}
	\phi_k(u) = \phi_\ell(u) = \sqrt{12} \left( u - \frac{1}{2} \right)
	\end{align}
for $u \in [0, 1]$ which results in $\rho_{k \ell}$ being 
the expected conditional 
Spearman correlation of $X$ and $Y$ given $Z$ with respect to the 
distribution of $Z$. Drawing inspiration from \eqref{eq:spearman} 
we define a class of functions $\phi : [0, 1] \to \R^q$ by letting
	\begin{align} \label{eq:trimmedspearmanphi}
	\phi_k(u) = c_k (u - m_k) \sigma_k(u)
	\end{align}
such that each $\phi_k : [0, 1] \to \R$ is determined by 
a Lipschitz continuous function $\sigma_k : [0, 1] \to \R$ 
with the support $\mathcal{T}_k$ of $\sigma_k$ a 
compact interval in $(0, 1)$, $\int_0^1 \sigma_k(u) \mathrm{d} u = 1$ and
	\begin{align*}
	m_k = \int u \sigma_k(u) \mathrm{d}u \quad \text{and} \quad c_k =   \left(\int (u - m_k)^2 \sigma_k(u)^2 \mathrm{d}u \right)^{-1/2}.
	\end{align*}
The choice \eqref{eq:trimmedspearmanphi} satisfies Assumption \ref{assump:phimap} 
(i) -- (iii) by construction, and if e.g. $\mathcal{T}_k \setminus \cup_{\ell \neq k} \mathcal{T}_{\ell} \neq \emptyset$ the functions are also linearly independent. We call the resulting generalized correlation 
$\rho$ a trimmed Spearman correlation, and we refer to the functions 
$\sigma_k$ as trimming functions. Note that if the supports $(\mathcal{T}_k)_{k=1}^q$ of 
$(\sigma_k)_{k=1}^q$ are choosen to be disjoint, then 
the covariance matrix $\Sigma$ of Theorem \ref{thm:rhoasympnormal} 
is the identity matrix.

A starting point for choosing a trimming 
function $\sigma$ is the normalized indicator
  \begin{align} \label{eq:indicator}
  u \mapsto (\lambda - \mu)^{-1} 1_{[\mu, \lambda]}(u)
  \end{align}
for $u \in [0, 1]$ where $0 < \mu < \lambda < 1$ are trimming 
parameters. However, \eqref{eq:indicator} is not a valid 
trimming function, since it is not Lipschitz. Therefore, 
we consider a simple linear approximation $\sigma : [0, 1] \to \R$ 
given by
  \begin{align} \label{eq:trimmingfunction}
  \sigma(u) = K f(u) \quad \text{and} \quad
  f(u) = \begin{cases}
  1, \ & u \in [\mu + \delta, \lambda - \delta] \\
  0, \ & u \in [\mu, \lambda]^c \\
  \delta^{-1} (u - \mu), \ & u \in [\mu, \mu + \delta)\\
  \delta^{-1} (\lambda - u), \ & u \in (\lambda - \delta, \lambda]
  \end{cases}
  \end{align}
and $K = (\lambda - \mu - \delta)^{-1}$. 
Here $0 < \delta < (\lambda - \mu) / 2$ is a fixed parameter that 
determines the accuracy of the approximation. It is elementary to verify that 
$\sigma$ given by \eqref{eq:trimmingfunction} is a valid trimming function, i.e., 
$\sigma$ is Lipschitz continuous 
with $\int \sigma(u) \d u = 1$ and support $[\mu, \lambda] 
\subset (0, 1)$.

The interpretation of a generalised correlation $\rho$ based on $\phi$ of the form 
\eqref{eq:trimmedspearmanphi} with trimming function $\sigma$ of the form 
\eqref{eq:trimmingfunction} is as follows. 
The entry $\rho_{k \ell}$ is an approximation of the expected conditional Spearman 
correlation between the observations of $X$ and $Y$, respectively, that lie in the 
$\mathcal{T}_k$-quantile range of the distribution of $X$ given $Z$ and the 
$\mathcal{T}_\ell$-quantile range of the distribution of $Y$ given $Z$, respectively, with 
respect to the distribution of $Z$. The matrix $\rho$ then summarizes this type of 
dependence within different quantile ranges of $X$ and $Y$ given $Z$.

\subsection{Practical considerations} \label{sec:practical}

Throughout Sections \ref{sec:correlationtest} and  \ref{sec:uniformlevelandpowerresults} we have 
analyzed our proposed test for conditional independence with an emphasis 
on modularity of the method regarding the choice of estimators $\hat{F}^{(n)}_{X \mid Z}$ 
and $\hat{F}^{(n)}_{Y \mid Z}$ of the conditional distribution functions and the choice of 
the function $\phi$ that defines the generalized correlation $\rho$ of Section \ref{sec:gcm}. 
This focus on the conceptual assumptions displays the generality of the method, but it also 
leaves the practitioner of conditional independence testing with a number of choices to be 
made. In this section we summarize a set of choices to make the method 
work out-of-the-box. 

Throughout the paper we have assumed that all random variables 
take values in the unit interval, i.e., $(X, Y, Z) \in [0,
1]^{d+2}$. This is not a restriction, since if e.g. $X \in \R$
we can always apply a strictly increasing, 
continuous transformation $ t: \R \to [0, 1]$ to obtain a new random variable 
$X' = t(X)$ with values in $[0, 1]$. Moreover, the initial conditional independence 
structure of $(X, Y, Z)$ is preserved since the transformation is marginal on $X$ and bijective.
The transformation $t$ can be chosen to be e.g. the logistic function.

In principle, an arbitrary and fixed marginal transformation
could be used for all variables, but we recommend to transform data to the unit
interval via marginal empirical distribution functions. This results in transformed 
variables known in the copula literature as pseudo copula observations. 
The transformation creates dependence, similar to the dependence created by
other common preprocessing techniques such as centering and scaling,
which the theoretical analysis has not accounted for. We suggest, nevertheless, 
to use this preprocessing technique in practise, and in the simulation study in Section 
\ref{sec:simulations} we use pseudo copula observations since it reflects how a practitioner would transform the variables.

To estimate the conditional distribution functions $\hat{F}^{(m, n)}_{X \mid Z}$ and $\hat{F}^{(m, n)}_{Y \mid Z}$ using Definition \ref{def:Fhat}, we suggest choosing
$\tau_{\min} = 0.01$ and $\tau_{\max} = 0.99$ and form the equidistant grid $(\tau_k)_{k=1}^m$ in $\mathcal{T} = [\tau_{\min}, \tau_{\max}]$ with 
the number of gridpoints $m = \lceil \sqrt{n} \rceil$. We then suggest 
using a model of the form
\eqref{eq:qrmodel} for both the quantile regression model $Q_{X \mid Z}(\tau_k \mid \cdot)$ and 
$Q_{Y \mid Z}(\tau_k \mid \cdot)$ for each $k=1, \dots, m$, where the bases $h_1$ and $h_2$ can be chosen to be e.g. an additive B-spline basis for each component of $Z$.

To test the hypothesis of conditional independence we suggest using the $\hat{\Psi}_n$ 
from Definition \ref{def:correlationtest} based on the estimated 
nonparametric residuals $(\hat{U}_{1, i}, \hat{U}_{2, i})_{i=1}^n$. 
To this end we choose $q \geq 1$ and let $\tau_{\min} = \lambda_0 < \cdots < \lambda_q = 
\tau_{\max}$ be an equidistant grid in $\mathcal{T}$. We then define the trimming function 
$\sigma_k$ to have the form \eqref{eq:trimmingfunction} with trimming parameters $\lambda_k$ and 
$\lambda_{k+1}$ and approximation parameter $\delta = 0.01 \cdot (\lambda_{k+1} - \lambda_k)$ 
for each $k=0, \dots, q-1$. We then define $(\sigma_k)_{k=1}^q$ according to \eqref{eq:trimmedspearmanphi}, 
compute the test statistic $\hat{\rho}_n$ using \eqref{eq:rhohat} and compute $\hat{\Psi}_n$ as in 
Definition \ref{def:correlationtest} using a desired significance level $\alpha \in (0, 1)$.

There are two non-trivial choices remaining. The first is the choice of 
bases $h_1$ and $h_2$ for the quantile regression models 
\mbox{$Q_{X \mid Z}(\tau_k \mid z) = h_1(z)^T \beta_{\tau_k}$} and 
\mbox{$Q_{Y \mid Z}(\tau_k \mid z) = h_2(z)^T \beta_{\tau_k}$}. Here the practitioner needs to make 
a qualified model selection. We recommend using a flexible basis such as an additive B-spline basis, and 
perform penalized estimation using \eqref{eq:betahat} to avoid overfitting. The second choice is the 
dimension of the generalized correlation $q \geq 1$, which corresponds to a choice of independence 
test in the partial copula. Note that the generalized correlation $\rho$ as 
above is defined for any $q \geq 1$, and there is 
conditional dependence, $X \notindep Y \mid Z$, if there exists $q \geq 1$ for which $\rho \neq 0$. 
We suggest trying one or a few, small values, e.g. $q \in \{1, \dots, 5\}$, and reject the hypothesis 
of conditional independence if one of the tests rejects the hypothesis, but of course be 
aware of multiple testing issues.

\section{Simulation Study} \label{sec:simulations}

In this section we examine the performance of our conditional independence test 
$\hat{\Psi}_n$ of Definition \ref{def:correlationtest}, when combining it with the 
quantile regression based conditional distribution function estimator $\hat{F}^{(m, n)}$ 
from Definition \ref{def:Fhat}. Firstly, we verify the level and power results obtained in 
Section \ref{sec:correlationtest} and Section \ref{sec:uniformlevelandpowerresults} empirically. 
Secondly, we compare the test with other conditional independence tests.
The test was implemented
in the R language \citep{R} using the \texttt{quantreg} package \citep{quantreg} as the
backend for performing quantile regression. The implementation and code for producing 
the simulations can be obtained from \url{https://github.com/lassepetersen/partial-copula-CI-test}. 

\subsection{Evaluation method}

We will evaluate the tests by their ability to hold level 
when data is generated from a distribution where conditional independence holds, and by 
their power when data is generated from a distribution where conditional 
independence does not hold. In order to make the results 
independent of a chosen significance level we will base the evaluation on 
the $p$-values of the tests. 

If a test has valid level, then we expect the $p$-value to be 
asymptotically $\mathcal{U}[0, 1]$-distributed. In Sections \ref{sec:level_and_power} 
and \ref{sec:comparisons} we evaluate the level by a Kolmogorov-Smirnov (KS) 
statistics as a function of sample size $n$, which is independent of any 
specific significance level. A small KS statistic is an indication of valid level.
To examine power we consider in Sections \ref{sec:level_and_power} 
and \ref{sec:comparisons} the 
$p$-values of the test as a function of the sample size, 
where we expect the $p$-values to 
tend to zero under the alternative of conditional dependence. 
Here a small $p$-value is an indication of large power.
In Section \ref{sec:local_alternative} we evaluate the power against a local 
alternative, which shrinks with the sample size $n$ toward the hypothesis 
of conditional independence with rate $n^{-\frac{1}{2}}$.

\subsection{Data generating processes} \label{sec:data_generation}

This section gives an overview of the data generating processes that we 
use for benchmarking and comparison. The first category consists of data generating processes of the form
  \begin{align} \tag{H} \label{eq:H}
  X = f_1(Z) + g_1(Z) \cdot \epsilon_{1}
  \quad \text{and} \quad 
  Y = f_2(Z) + g_2(Z) \cdot \epsilon_{2}
  \end{align}
where $f_1, f_2, g_1, g_2 : \R^d \to \R$ belong to some function class and 
$\epsilon_1, \epsilon_2$ are independent errors. For data generating processes of 
type \eqref{eq:H}, conditional independence is satisfied.
The second category consists of data generating processes of the form 
  \begin{align} \tag{A} \label{eq:A}
  X = f_1(Z) + g_1(Z) \cdot \epsilon_{1}
  \quad \text{and} \quad 
  Y = f_2(Z, X) + g_2(Z, X) \cdot \epsilon_{2}
  \end{align}
where again $f_1, g_1 : \R^{d} \to \R$ and $f_2, g_2 : \R^{d+1} \to \R$ 
belong to some function class and 
$\epsilon_1, \epsilon_2$ are independent errors. Under data generating processes
of type \eqref{eq:A}, conditional independence is not satisfied. We will 
consider four different data generating processes corresponding to different choices 
of functions $f_1$, $g_1$, $f_2$ and $g_2$ and error distributions.

\begin{itemize}
\item[(1)] For data generating processes $\mathrm{H}_1$ and $\mathrm{A}_1$ we let
  \begin{align*}
  f_k(w_1, \dots, w_d) & = \sum_{j=1}^d
  \beta_{1, k, j} w_j + \beta_{2, k, j} w_{j}^2
  \\
  g_k(w_1, \dots, w_d) & = \exp \left(
  - \left|
  \sum_{j=1}^d
  \alpha_{1, k, j} w_j + \alpha_{2, k, j} w_{j}^2
  \right|
  \right)
  \end{align*}
for $k=1, 2$ and real valued coefficients 
$(\alpha_{\ell, k, j}, \beta_{\ell, k, j})_{\ell=1, 2, k=1,2, j=1, \dots, d}$.
Here each $Z_{j} \sim \mathcal{U}[-1, 1]$ 
independently, $\epsilon_{1}$ follows an asymmetric Laplace distribution with location 
$0$, scale $1$ and skewness $0.8$, and $\epsilon_{2}$ follows
a Gumpel distribution with location $0$ and scale $1$.

\item[(2)] For data generating processes $\mathrm{H}_2$ and $\mathrm{A}_2$ we let 
$g_1 = g_2 = 1$ and 
  \begin{align*}
  f_k(w_1, \dots, w_d) = \sum_{j=1}^d \beta_{k, j} w_j
  \end{align*}
for $k=1, 2$ and real valued coefficients $(\beta_{k, j})_{k=1, 2, j=1, \dots, d}$. 
Here each $Z_{j} \sim \mathcal{U}[-1, 1]$ 
independently and both $\epsilon_1$ and $\epsilon_2$ follow a $\mathcal{N}(0, 1)$-distribution 
independently.

\item[(3)] For data generating processes $\mathrm{H}_3$ and $\mathrm{A}_3$ we let 
$g_1 = g_2 = 1$ and 
  \begin{align*}
  f_k(w_1, \dots, w_d) = \sum_{j=1}^d \beta_{1, k, j} w_j + \beta_{2, k, j} w_j^2
  \end{align*}
for $k=1, 2$ and real valued coefficients 
$(\beta_{\ell, k, j})_{\ell=1, 2, k=1, 2, j=1, \dots, d}$. 
Here each $Z_{j} \sim \mathcal{U}[-1, 1]$ 
independently and both $\epsilon_1$ and $\epsilon_2$ follow a $\mathcal{N}(0, 1)$-distribution 
independently.

\item[(4)] For data generating processes $\mathrm{H}_4$ and $\mathrm{A}_4$ we let 
$f_1 = f_2 = 0$ and 
  \begin{align*}
  g_k(w_1, \dots, w_d) = \sum_{j=1}^d \beta_{1, k, j} w_j + \beta_{2, k, j} w_j^2
  \end{align*}
for $k=1, 2$ for real valued coefficients 
$(\beta_{\ell, k, j})_{\ell=1, 2, k=1, 2, j=1, \dots, d}$. 
Here each $Z_{j} \sim \mathcal{U}[-1, 1]$ 
independently and both $\epsilon_1$ and $\epsilon_2$ follow a $\mathcal{N}(0, 1)$-distribution 
independently.
\end{itemize}

Each time we simulate from data generating processes $\mathrm{H}_1, \dots, \mathrm{H}_4$ 
we first draw the coefficients of the functions $f_k, g_k$ from a $\mathcal{N}(0, 1)$-distribution 
in order to make the results independent of a certain combination of parameters. When we 
simulate from the data generating processes $\mathrm{A}_1, \mathrm{A}_2$ and $
\mathrm{A}_3$ we first draw the 
coefficients of $f_k, g_k$ to be either $-1$ or $1$ with equal probability in order to 
fix the signal to noise ratio between the predictors and responses. When simulating from 
$\mathrm{A}_4$ we simulate the coefficients of $g_k$ to be either $-5$ or $5$, because the 
conditional dependence lies in the variance for $\mathrm{A}_4$, and a stronger signal is 
needed to compare the power of the tests using the same samples sizes as for 
$\mathrm{A}_1, \mathrm{A}_2$ and $\mathrm{A}_3$.

The data generating processes $\mathrm{H}_2, \mathrm{H}_3, \mathrm{H}_4, 
\mathrm{A}_2, \mathrm{A}_3$ and $\mathrm{A}_4$ can be shown to satisfy Assumption 
\ref{assump:highdimquantreg} that is needed for 
Corollary \ref{cor:Qhatconsistency}, since they are linear (in the parameters) 
location-scale models with bounded covariates \citep[Section 2.5]{belloni2011}. The processes
$\mathrm{H}_1$ and $\mathrm{A}_1$ are not of this form, since $g_1$ and $g_2$ are 
nonlinear in the parameters. However, we include these in the simulation study 
to test the robustness of the test.

\subsection{Level and power of partial copula test} \label{sec:level_and_power}

\begin{figure}[t]
\centering
\includegraphics[scale=0.6]{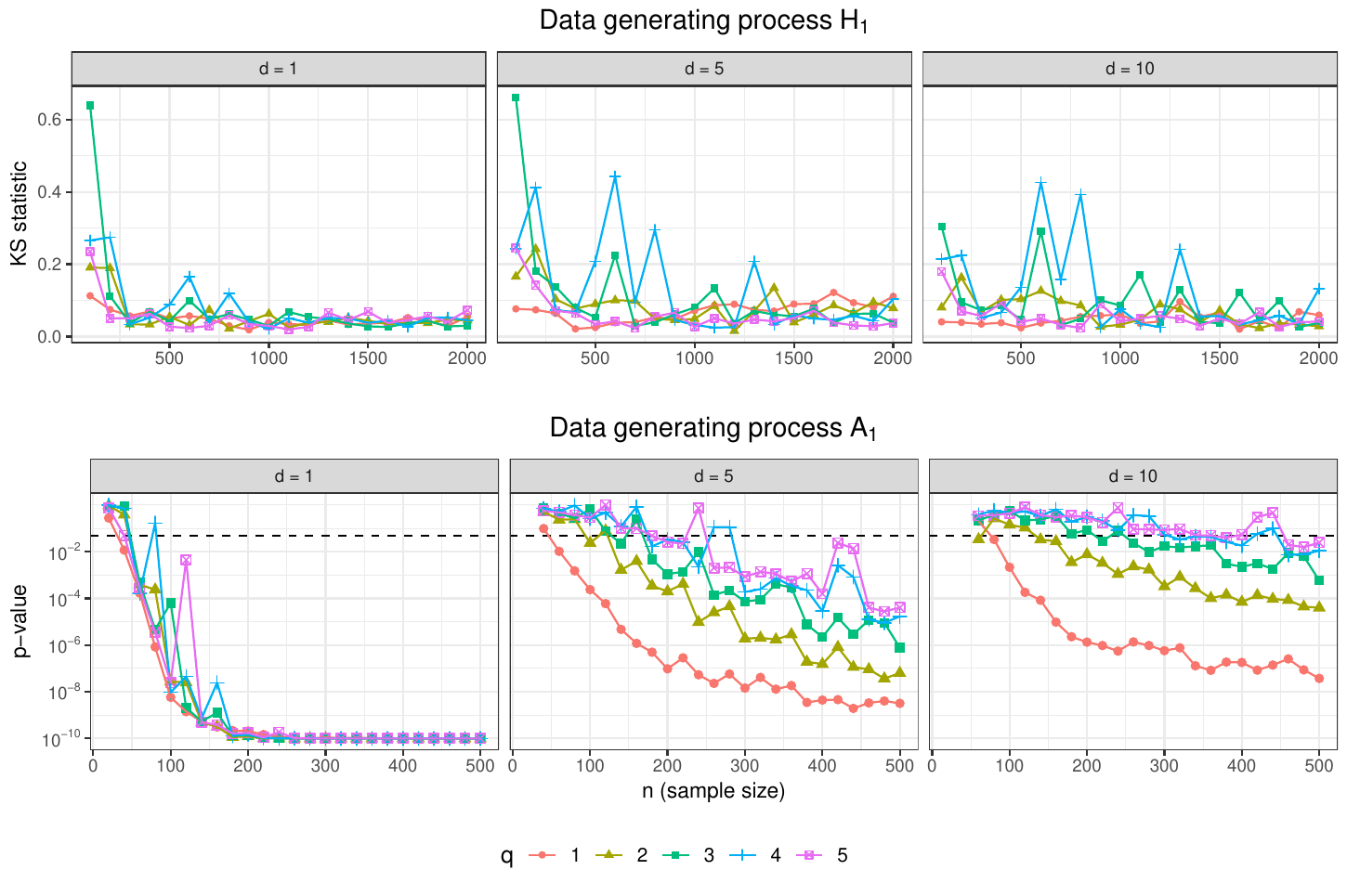} 
\caption{Top: KS statistic for equality with a $\mathcal{U}[0, 1]$-distribution 
of the p-values of the two tests 
computed from 500 simulations from $\mathrm{H}_1$ for each 
combination of $n$ and $d$. Bottom: Average p-values of the 
two tests separately computed over 200 simulations from $\mathrm{A}_1$ for 
each combination of $n$ and $d$. Dashed line indicates the common significance level 
$0.05$. For visual purposes all $p$-values have been truncated at $10^{-10}$.}
\label{fig:plot1}
\end{figure}

In this section we examine the level and power properties of 
the test $\hat{\Psi}_n$. 
We examine the performance of the test on data 
generating processes $\mathrm{H}_1$ and $\mathrm{A}_1$ for 
dimensions $d \in \{1, 5, 10\}$ of $Z$. 
The test is performed as described 
in Section \ref{sec:practical}. 
As the quantile regression model we 
use an additive model with a B-spline basis 
of each variable with 5 degrees of 
freedom, and we try $q \in \{1, \dots, 5\}$.
The result of the simulations 
can be seen in Figure \ref{fig:plot1}. 
We observe that for $d=1$ all five tests obtain level asymptotically 
under $\text{H}_1$, 
while for higher dimension $d \in \{5, 10\}$ the test with 
$q = 4$ has minor problems holding level. We also see that the 
$p$-values for all five tests tend to zero as the sample size increases 
under $\text{A}_1$. The convergence rate of the $p$-value depends on 
the dimension $d$ such that a higher dimension gives a slower convergence 
rate. In conclusion we observe that our test holds level under 
a complicated data generating distribution ($\text{H}_1$), where there 
is a nonlinear conditional mean and variance dependence and skewed error 
distributions with super-Gaussian tails. Moreover, the test has power 
against the alternative of conditional dependence ($\text{A}_1$), however, 
for $d=1$ we see that $q = 5$ gives the best power, while $q=1$ 
gives the best power for $d \in \{5, 10\}$. The testing procedure also 
displays robustness to the fact that the quantile regression models are 
misspecified.

\subsection{Comparison with other tests} \label{sec:comparisons}

We now compare the partial copula based test $\hat{\Psi}_n$ with other 
nonparametric tests. 
We will compare with a residual based method, since 
this is another class of conditional independence test based on 
nonparametric regression. In order to describe this test we let
  \begin{align*}
  R_{1, i} = X_i - \hat{f}(Z_i) \ \quad \text{and} 
  \quad 
  R_{2, i} = Y_i - \hat{g}(Z_i)
  \end{align*}
for $i=1, \dots, n$ be the residuals obtained 
when performing conditional mean regression $\hat{f}$ of 
$f(z) = E(X \mid Z=z)$ and $\hat{g}$ of 
$g(z) = E(Y \mid Z=z)$ obtained from a 
sample $(X_i, Y_i, Z_i)_{i=1}^n$.
We compare the following conditional independence tests:

\begin{itemize}
\item \textbf{GCM}: The Generalised Covariance Measure  
which tests for vanishing correlation between the 
residuals $R_1$ and $R_2$ given as above \citep{shah2018hardness}. 

\item \textbf{NPN correlation}: Testing for vanishing partial correlation 
in a nonparanormal distribution \citep{harris2013pc}. This is a generalization of the 
partial correlation, which assumes a Gaussian dependence structure, but allows for arbitrary 
marginal distributions.

\item \textbf{PC}: Our partial copula based test $\hat{\Psi}_n$ for 
$q \in \{1, 3, 5\}$ as described in Section \ref{sec:practical}.
\end{itemize}
\begin{figure}[t]
\centering
\includegraphics[scale=0.6]{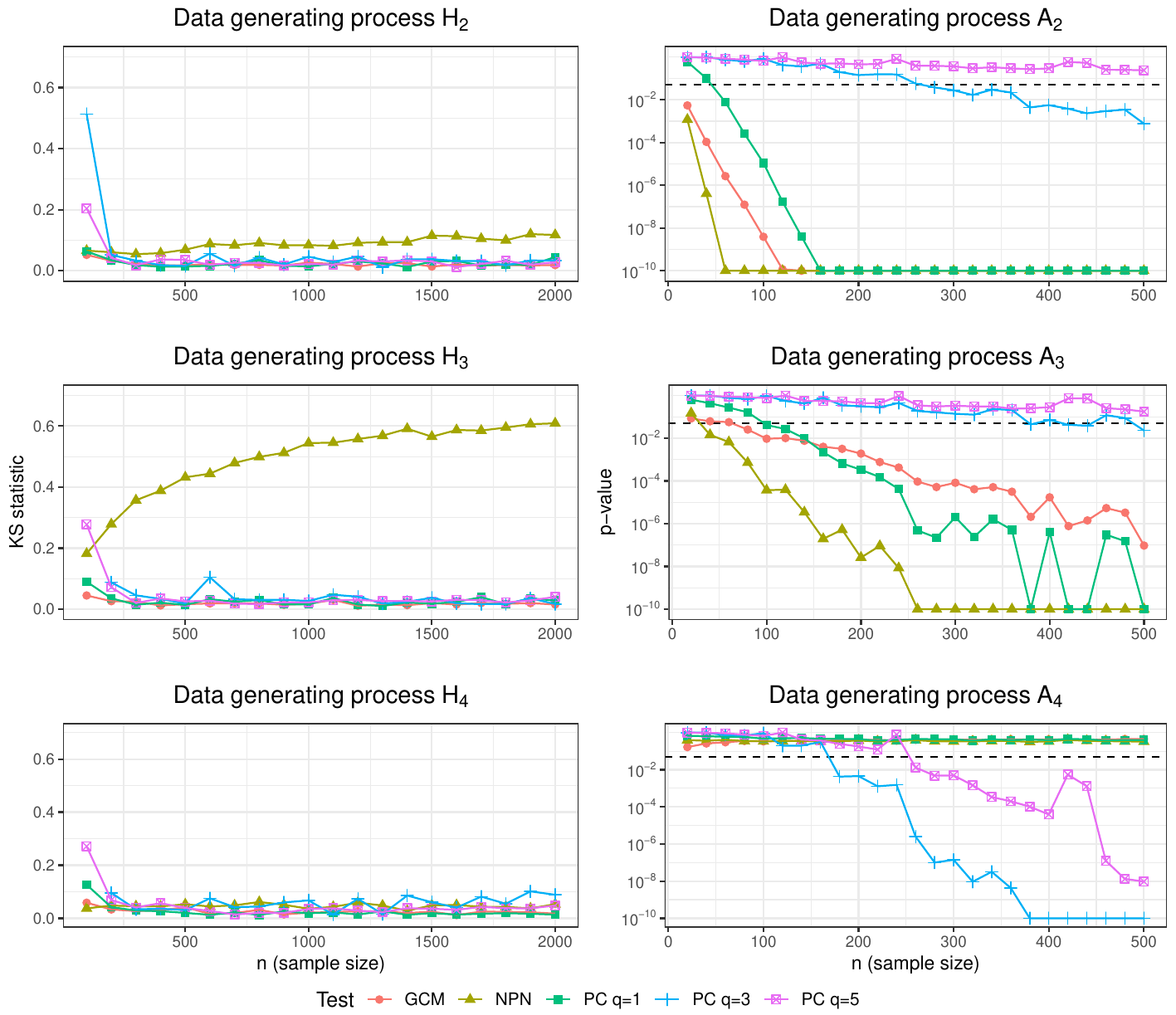} 
\caption{Left column: KS statistic for equality with a
$\mathcal{U}[0, 1]$-distribution of the p-values of the five tests computed from 
500 simulations from $\mathrm{H}_2, \mathrm{H}_3$ and $\mathrm{H}_4$, respectively, for 
each sample size $n$. Right column:
Average p-values of the five tests separately computed over 200 simulations from 
$\mathrm{A}_2, \mathrm{A}_3$ and $\mathrm{A}_4$, respectively, for each sample size 
$n$. For all simulations the dimension is fixed at $d=3$.
Dashed line indicates the common significance level 
$0.05$. For visual purposes all $p$-values have been truncated at $10^{-10}$.}
\label{fig:plot2}
\end{figure}

We consider the behavior of the tests under 
$\mathrm{H}_2, \mathrm{A}_2, \mathrm{H}_3, \mathrm{A}_3, \mathrm{H}_4$ and $\mathrm{A}_4$.
For fairness of comparison we choose  our quantile and mean regression models 
to be the correct model class such that the tests perform at their oracle level, 
e.g., for $\mathrm{H}_3$ 
we fit additive models with polynomial basis of degree $2$. We fix the dimension $d$ of $Z$ 
to be $3$ in all simulations for simplicity. 
The results of the simulations can be seen in Figure \ref{fig:plot2}. 

\begin{figure}[t]
\centering
\includegraphics[scale=0.51]{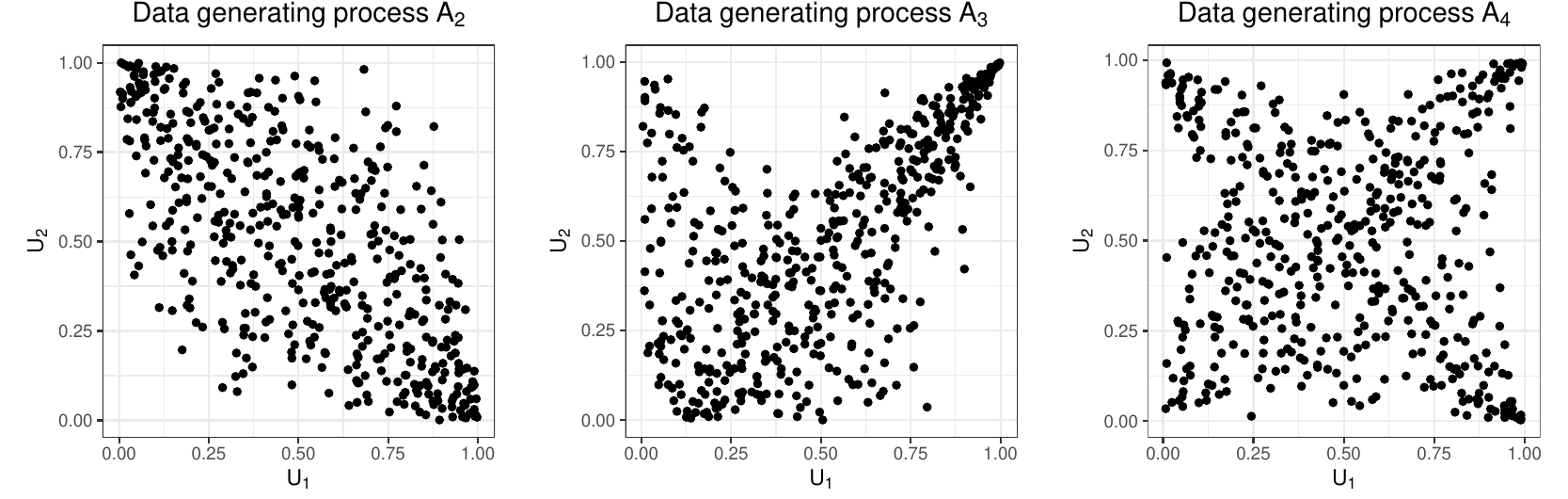}
\caption{Estimated partial copula 
$(\hat{U}_{1, i}, \hat{U}_{2, i})_{i=1}^n$ from one realization 
from each of the data generating processes $\text{A}_2, \text{A}_3$ and 
$\text{A}_4$ for $d = 3$ and $n=500$.}
\label{fig:plot3}
\end{figure}

Under $\text{H}_2$ all five tests hold level, and 
we see that the NPN test has greatest power against $\text{A}_2$ followed 
by the GCM and $\hat{\Psi}_n$ with $q=1$, while $\hat{\Psi}_n$ with 
$q \in \{3, 5\}$ does not have much power against $\text{A}_2$. In order to 
intuitively understand the effect of $q$ 
see Figure \ref{fig:plot3}. We see that in the estimated 
partial copula the dependence is captured by the overall correlation, 
while dividing 
$[\tau_{\min}, \tau_{\max}] \times [\tau_{\min}, \tau_{\max}]$ into 
subregions does not reveal finer dependence structure. Hence 
$q = 1$ is suitable to detect the dependence for $\text{A}_1$. 

Under $\text{H}_3$ both the GCM test and $\hat{\Psi}_n$ with 
$q \in \{1, 3, 5\}$ hold level, but the 
NPN test does not hold level under $\text{H}_3$, which is due to 
the nonlinear response-predictor relationship. However, since both 
the GCM and $\hat{\Psi}_n$ test takes the nonlinearity into account, they 
can effectively filter away the $Z$-dependence. The NPN test has 
greatest power against the alternative $\text{A}_3$
following by $\hat{\Psi}_n$ with $q=1$ and the GCM test. In 
Figure \ref{fig:plot3} we again see that the dependence in the 
estimated partial copula is described by the overall correlation, 
while dividing into subregions results a generalized correlation 
with elements that are close to zero, i.e., here 
$q=1$ is suitable for capturing the dependence. 

Under $\text{H}_4$, all test hold level. Note that the NPN test holds level even 
though there is a nonlinear conditional variance relation, since this
is still a nonparanormal distribution. 
We also see that neither the GCM test nor the NPN test has 
power against $\text{A}_4$, while $\hat{\Psi}_n$ has some power against 
$\text{A}_4$ with the greatest power for $q = 3$. In Figure \ref{fig:plot3} 
we see that there is a clear dependence in the estimated partial copula, 
but that the overall correlation is close to zero. However, when 
dividing into subregions the generalized correlation is able 
to detect the dependencies in the tails of the distributions. 
\subsection{Power under local alternatives} \label{sec:local_alternative}

Though GCM did not have power against the specific alternative $\text{A}_4$,
it maintains level and it has power against a broad class of alternatives. 
To understand better when $\hat{\Psi}_n$ can be expected to have greater power than GCM, we 
consider a simulation, which is a small variation of the simulations 
presented in Section \ref{sec:data_generation}. 

The dimension is fixed as $d = 1$, $Z \sim \mathcal{U}([0,1])$ is uniformly 
distributed on $[0, 1]$, $\epsilon_{1}$, $\epsilon_{2}$ and $W$ are independent
and $\mathcal{N}(0, 1)$-distributed, and 
\begin{align} \tag{A} \label{eq:A}
X = (\beta Z^2 + 1) \epsilon_{1} + \gamma W
\quad \text{and} \quad 
Y = (\beta Z^2 + 1) \epsilon_{2} + \gamma W
\end{align}
for parameters $\beta, \gamma \in \mathbb{R}$. Conditionally on $Z$, the distribution 
of $(X, Y)$ is a bivariate Gaussian distribution, and $X$ and $Y$ are 
conditionally independent if and only if $\gamma = 0$. 
\begin{figure}[t]
\centering
\includegraphics[scale=0.55]{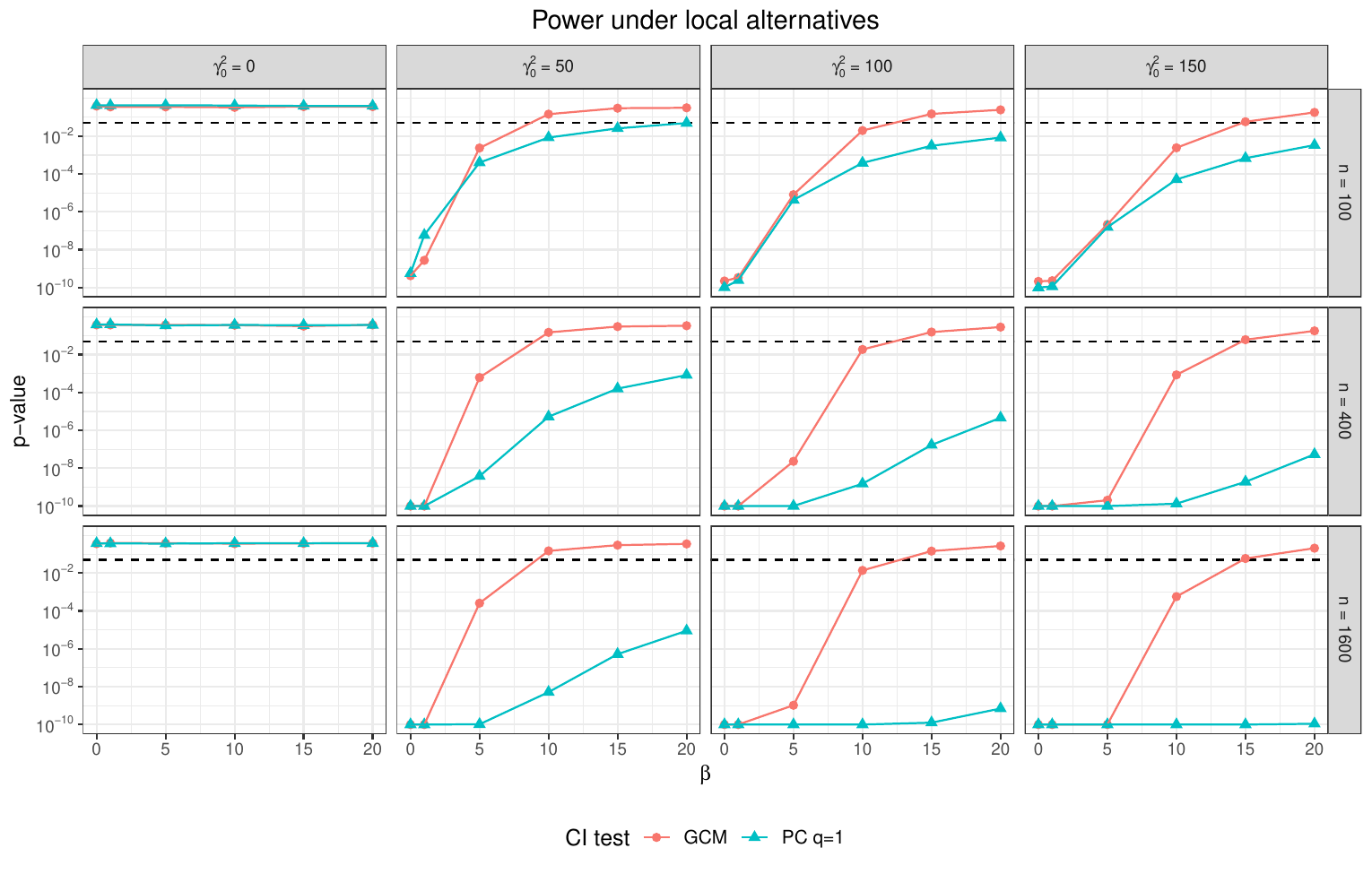} 
\caption{Average p-values of the three tests separately computed over 500 
simulations for different values of the parameters $\beta$ and $\gamma_0$ and 
different sample sizes $n$. The value $\gamma_0^2 = 0$ is equivalent to 
conditional independence. Other values correspond to the local alternatives
$\gamma^2 = \gamma_0^2 / \sqrt{n}$. The values of $\beta$ determine the 
amount of heterogeneity of the mean regression residual variances with 
$\beta = 0$ meaning constant residual variance and $\beta = 20$ meaning 
substantial heterogeneity.
Dashed line indicates the common significance level 
$0.05$. For visual purposes all $p$-values have been truncated at $10^{-10}$.}
\label{fig:plot4}
\end{figure}
We examine level and power by simulating 500 data sets for sample 
sizes $n \in \{100, 400, 1600\}$ and all combinations of parameters 
$\beta \in \{0, 1, 5, 10, 15, 20\}$, and local alternatives
$$\gamma^2 = \frac{\gamma_0^2}{\sqrt{n}}$$
for $\gamma_0^2 \in \{0, 50, 100, 150\}$. Note that 
$f(z) = E(X \mid Z = z) = 0$ and $g(z) =  E(Y \mid Z = z) = 0$, 
which is exploited for GCM instead of 
estimating $f$ and $g$. This should only increase the power of 
GCM relative to fitting any model of the conditional expectations. 
We perform the test $\hat{\Psi}_n$ as described in Section 
\ref{sec:practical} using $q=1$, and the quantile regression model 
is fitted using a polynomial basis of degree 2.

Figure \ref{fig:plot4} shows the results of the simulation. 
Both GCM and $\hat{\Psi}_n$ maintain level for $\gamma_0^2 = 0$. 
$\hat{\Psi}_n$ has comparable 
or superior power relative to GCM in all other cases. 
Both tests have decreasing power as a function of 
$\beta$, but $\hat{\Psi}_n$ maintains power even for large values of $\beta$, where GCM 
has almost no power. The power of $\hat{\Psi}_n$ against the 
local alternatives increases 
with the sample size, which shows how the increased precision for larger samples
of the quantile regression based distribution functions improves power. We do not see
the same for GCM, partly because no mean value model is fitted. 

As $\beta$ quantifies the conditional variance heterogeneity of $X$ and $Y$ 
given $Z$, we conclude that though GCM remains a valid test under conditional 
variance heterogeneity, its test statistic does not adequately account for 
the heterogeneity, and GCM has inferior power 
under local alternatives when compared to $\hat{\Psi}_n$. 

\section{Discussion} \label{sec:discussion}

The first main contribution of this paper is an estimator 
of conditional distribution functions $\hat{F}^{(m, n)}$ based on 
quantile regression. We have shown that the estimator is pointwise 
(uniformly) consistent over a set of distributions $\mathcal{P}_0 \subset 
\mathcal{P}$ given that the quantile regression procedure is 
pointwise (uniformly) consistent over $\mathcal{P}_0$. 
Moreover, we showed that the convergence 
rate of the quantile regression procedure 
can be transferred directly to the estimator $\hat{F}^{(m, n)}$.

The second main contribution of this paper is an analysis of a 
nonparametric test for conditional independence based on the 
partial copula construction. We introduced a class of tests given in 
terms of a generalized correlation dependence measure $\rho$ with the 
leading example being a trimmed version of the Spearman correlation. 
We showed that the test achieves asymptotic pointwise (uniform) level 
and power over $\mathcal{P}_0$ given that the conditional distribution function estimators 
are pointwise (uniformly) consistent over $\mathcal{P}_0$ 
with rate functions 
$g_P$ and $h_P$ satisfying $\sqrt{n} g_P(n) h_P(n) \to 0$. 
The partial copula has previously been considered for 
conditional independence testing in the literature, however, to 
the best of our knowledge, the results presented here are the first to 
explicitly connect the consistency requirements of the conditional 
distribution function estimators to level and power properties of the 
test.

Lastly, we established through a simulation study that the proposed  test is
sound under complicated data generating distributions, and that  it has power
comparable to or even better than other  state-of-the-art nonparametric
conditional independence tests.  In particular, we demonstrated that our test
has superior power against alternatives with variance heterogeneity between $X$
and $Y$ given $Z$ when compared to conditional independence tests based on
conventional residuals. We note that due to Daudin's lemma, tests based on
conventional residuals can obtain power against any alternative if suitable 
transformations of $X$ and $Y$ are considered. In particular, if 
$X^2$ and $Y^2$ were used in our simulation study, GCM would have power against 
$\text{A}_4$. We tested the use of GCM in combination with $X^2$ and $Y^2$ 
in all our simulations (data not shown), and though it had some power against
$\text{A}_4$, it was comparable to or inferior to just using GCM in all 
other simulations. Thus to obtain good power properties, the specific choice 
of transformation appears important and to depend on the data generating 
distribution. 

An important point about the test is the rate requirement 
$\sqrt{n} g_P(n) h_P(n) \to 0$ needed to achieve asymptotic level. 
The product structure means that the test is sound under 
quantile regression models with slower consistency rates than the 
usual parametric $n^{-1/2}$-rate. This opens up the methodology 
to nonparametric machine learning models. An interesting direction 
of research would be to empirically assess the performance of the 
test using machine learning inspired 
quantile regression models, such as deep neural networks, where explicit 
consistency rates are not available. We hypothesize 
that the method will still perform well in these scenarios due to the 
weak consistency requirement.

In this paper we have considered univariate $X$ and $Y$. A possible extension of 
the test is to allow $X \in [0, 1]^{r_1}$ and $Y \in [0, 1]^{r_2}$ with 
$r_1, r_2 \geq 1$, and then consider the nonparametric residual $U_1 \in [0, 1]^{r_1}$ 
of $X$ given $Z$ by performing coordinatewise probability integral transformations 
$U_{1, k} = F_{X_k \mid Z}(X_k \mid Z)$ for $k = 1, \dots, r_1$, and similarly for 
constructing the nonparametric residual $U_2 \in [0, 1]^{r_2}$ of 
$Y$ given $Z$. Conditional independence 
$X \indep Y \mid Z$ then implies pairwise independence of $U_{1,k}$ and $U_{2,l}$
for each $k = 1, \ldots, r_1$ and $l = 1, \ldots, r_2$. Combining our proposed 
test statistics for each such pair yields an $r_1r_2 q^2$-dimensional 
test statistic, whose distribution under the hypothesis of conditional independence
will be asymptotically Gaussian with mean 0. Its covariance matrix will
only be partially known, though, due to the potential dependence between the pairs, 
but the unknown part could be estimated from the estimated nonparametric residuals. 
The multivariate statistic could be aggregated into a univariate test statistic 
in various ways, e.g. by a quadratic transformation as in \eqref{eq:chistatistic}, 
or by the maximum of the absolute values of its coordinates. In the low-dimensional
case for fixed $r_1$ and $r_2$ our results would carry over immediately, and 
we expect that using the maximum could lead to high-dimensional results similar to 
Theorem 9 by \cite{shah2018hardness}. 

A key property of the partial copula is that the nonparametric residuals 
$U_1$ and $U_2$ are independent under conditional independence and not only uncorrelated, 
which is the case for conventional residuals in additive noise models. Therefore, an 
important question is whether asymptotic level and power 
guaranties can be proven, when combining the partial copula with 
more general independence tests. In this paper we have focused 
on dependence measures of the form $\rho = E_P(\phi(U_1) \phi(U_2)^T)$ 
and tests based on 
	\begin{align*}
	\hat{\rho}_n = \frac{1}{n} \sum_{i=1}^n \phi(\hat{U}_{1, i}) 
	\phi(\hat{U}_{2, i})^T
	\end{align*}
because it gives a flexible and general test for independence in the 
partial copula, it can be computed in linear time in the size of data, and most importantly 
its asymptotic theory is standard and easy to establish and apply. It also 
clearly illustrates the 
transfer of consistency of the conditional distribution function estimators 
to properties of the test. It is ongoing work to establish a parallel asymptotic 
theory for dependence measures of the form 
$\theta = E_P(h(U_1, U_2))$, where $h$ is a kernel 
function, and whose estimators are $U$-statistics. This could potentially 
yield more powerfull tests against complicated alternatives of conditional 
dependence, but at the prize of increased computational complexity.


\acks{This work was supported by a research grant (13358) from VILLUM FONDEN.}



\appendix

\section{Proofs} \label{appendix:proofs}

This appendix gives proofs of the main results of the paper. 
Throughout the proofs we will ignore the dependence of certain 
terms on the sample size to ease notation, e.g. we write 
$\hat{U}_{1, i}$ instead of $\hat{U}^{(n)}_{1, i}$ and $\hat{q}_{k, z}$ instead of 
$\hat{q}_{k, z}^{(n)}$. 

\subsection{Proof of Proposition \ref{prop:Ftildeapproximation}}

We need to bound the supremum 
	\begin{align*}
	||F - \tilde{F}^{(m)}||_{\mathcal{T}, \infty} = \sup_{z \in [0, 1]^d} 
	\sup_{t \in Q(\mathcal{T} \mid z)} |F(t \mid z) - \tilde{F}^{(m)}(t \mid z) |.
	\end{align*}
First we fix $z \in [0, 1]^d$ and inspect the inner supremum. 
By construction we have 
  \begin{align*}
  F(q_{k, z} \mid z) = \tilde{F}^{(m)}(q_{k, z} \mid z) = \tau_k
  \end{align*}
for $k = 1, \dots, m$. Furthermore, since both $F$ and $\tilde{F}^{(m)}$ 
are continuous and increasing in $t \in [0, 1]$ we have that 
  \begin{align*}
  \sup_{t \in [q_{k, z}, q_{k+1, z}]} | F(t \mid z) - \tilde{F}^{(m)}(t \mid z) |
  \leq \tau_{k+1} - \tau_k
  \end{align*}
for each $k = 1, \dots, m-1$. Since $Q(\mathcal{T} \mid z) = [q_{\min, z}, q_{\max, z}] = 
\bigcup_{k=1}^{m-1} [q_{k, z}, q_{k+1, z}]$ 
we now have 
  \begin{align*}
  \sup_{t \in Q(\mathcal{T} \mid z)} 
  | F(t \mid z) - \tilde{F}^{(m)}(t \mid z)| 
  & = \max_{k = 1, \dots, m-1} 
  \sup_{t \in [q_{k, z}, q_{k+1, z}]} | F(t \mid z) - \tilde{F}^{(m)}(t \mid z) |
  \\
  & \leq 
  \max_{k=1, \dots, m-1} (\tau_{k+1} - \tau_k) = \kappa_m.
  \end{align*}
The result now follows from taking supremum over $z \in [0, 1]^d$ as the 
right hand side of the inequality does not depend on $z$. \hfill $\square$

\subsection{Proof of Proposition \ref{prop:Fhatconsistent}}

We need to bound the supremum
  \begin{align*}
  \|\tilde{F}^{(m)} - \hat{F}^{(m, n)}\|_{\mathcal{T}, \infty} 
  = \sup_{z \in [0, 1]^d} \sup_{t \in Q(\mathcal{T} \mid z)}
  |\tilde{F}^{(m)}(t \mid z) - \hat{F}^{(m, n)}(t \mid z)|.
  \end{align*}
Our proof strategy is the following. First we evaluate the inner supremum over 
$t \in Q(\mathcal{T} \mid z)$ analytically to obtain a bound in terms of the 
quantile regression prediction error. Then we will evaluate the outer supremum over 
$z \in [0, 1]^d$ and use the assumed consistency from 
Assumption \ref{assump:Qhatconsistency}. First define the two quantities 
  \begin{align*}
  A(m, n, z) := \kappa_m \cdot 
  \frac{\max_{k=1, \dots, m} |q_{k, z} - \hat{q}^{(n)}_{k, z}|}
  {\min_{k=1, \dots, m-1}(q_{k+1, z} - q_{k, z})}
  \end{align*}
and
  \begin{align*}
  B(m, n, z) := \kappa_m \cdot 
  \frac{\max_{k=1, \dots, m} |q_{k, z} - \hat{q}^{(n)}_{k, z}|}
  {\min_{k=1, \dots, m-1}(\hat{q}^{(n)}_{k+1, z} - \hat{q}^{(n)}_{k, z})}.
  \end{align*}
We then have the following key result regarding the inner supremum over $t \in 
Q(\mathcal{T} \mid z)$. 

\begin{proposition} \label{prop:innersupremum}

Let Assumption \ref{assump:Qhatconsistency} (i) be satisfied. Then for 
all $P \in \mathcal{P}_0$ and $\epsilon > 0$ there exists $N \geq 1$ such 
that for all $n \geq N$,
  \begin{align*}
  \sup_{t \in Q(\mathcal{T} \mid z)} |\tilde{F}^{(m)}(t \mid z) - \hat{F}^{(m, n)}(t \mid z)|
  \leq \max \{ A(m, n, z), B(m, n, z) \}
  \end{align*}
for all $z \in [0, 1]^d$ and all grids $(\tau_k)_{k=1}^m$ in $\mathcal{T}$
with probability at least $1 - \epsilon$. 
\end{proposition}

We need a number of auxilliary results 
before proving Proposition \ref{prop:innersupremum}. 
We start by proving the following key lemma that reduces the number of 
distinct cases of relative positions of the true conditional quantiles 
$q_{k, z}$ and the estimated conditional quantiles $\hat{q}_{k, z}$.

\begin{lemma} \label{lemma:distinctquantilepositions}
Let Assumption \ref{assump:Qhatconsistency} (i) 
be satisfied. Then for each $P \in \mathcal{P}_0$ and $\epsilon > 0$ there
exists $N \geq 1$ such that for all $n \geq N$ we have that 
$\hat{q}_{k, z} \in (q_{k-1, z}, q_{k+1, z})$ for each $k=1, \dots, m$ and 
$z \in [0, 1]^d$ and for all grids $(\tau_k)_{k=1}^m$ in $\mathcal{T}$ 
with probability at least $1 - \epsilon$.
\end{lemma}

\begin{proof}
Fix a distribution $P \in \mathcal{P}_0$. 
Let $\mathcal{G}$ be the set of all grids $(\tau_k)_{k=1}^m$ in $\mathcal{T}$. 
Then
  \begin{align*}
  \sup_{\mathcal{G}} \sup_{z \in [0, 1]^d} \max_{k=1, \dots, m}
  |q_{k, z} - \hat{q}_{k, z}^{(n)}| \leq 
  \sup_{z \in [0, 1]^d} \sup_{\tau \in \mathcal{T}} 
  |Q(\tau \mid z) - \hat{Q}^{(n)}(\tau \mid z)| \stackrel{P}{\to} 0
  \end{align*}
under Assumption \ref{assump:Qhatconsistency} (i). Since 
$q_{k, z} \in (q_{k-1, z}, q_{k+1, z})$ for each $k=1, \dots, m$ and $z \in [0, 1]^d$
for all grids $(\tau_k)_{k=1}^m$ in $\mathcal{T}$ the result follows. 
\end{proof}

Next we have some lemmas giving the supremum of certain functions over certain intervals that will 
be useful in the main proof.

\begin{lemma} \label{lemma:basicinequality1}
Let $a \leq b < c \leq d$ and $f(t) = \frac{t - a}{c - a} - \frac{t - b}{d - b}$. 
Then $\sup_{t \in [b, c]} f(t) = \max\{ \frac{b-a}{c-a}, \frac{d-c}{d-b} \}$.
\end{lemma}

\begin{proof}
Note that $f$ is a linear function. Thus the supremum is
obtained in one of the intervals endpoints, i.e.,
$\sup_{t \in [b, c]} f(t) = \max\{f(b), f(c)\}$. We see that
  \begin{align*}
  f(b) = \frac{b-a}{c-a}
  \quad \text{and} \quad
  f(c) = 1 - \frac{c-b}{d-b}
  =
  1 - \frac{c - d + d - b}{d - b}
  =
  \frac{d - c}{d - b},
  \end{align*}
which shows the result. 
\end{proof}

\begin{lemma} \label{lemma:basicinequality2}
Let $a < b \leq c < d$ and
$f(t) = \alpha + \beta \cdot \frac{t-b}{d-b} - \alpha \cdot \frac{t-a}{c-a}$ where $\alpha, \beta > 0$.
Then we have
$\sup_{t \in [b, c]} f(t)= \max \{ \alpha \cdot \frac{c-b}{c-a}, \beta \cdot \frac{c-b}{d-b}\}$.
\end{lemma}

\begin{proof}
The function $f$ is a linear function, and hence the supremum is 
obtained in one of the interval endpoints. We see that
  \begin{align*}
  f(b) = \alpha - \alpha \cdot \frac{b-a}{c-a} = \alpha \cdot \frac{c-b}{c-a} 
  \quad \text{and} \quad
  f(c) = \beta \cdot \frac{c-b}{d-b},
  \end{align*}
which shows the claim.
\end{proof}

\begin{lemma} \label{lemma:basicinequality3}
Let $a \leq b < c < d$ and $f(t) = |g(t)|$ where 
$g(t) = \frac{t-b}{c-b} - \frac{t - a}{d - a}$. Then we have that
$\sup_{t \in [b, c]} f(t) = \max\{\frac{b - a}{d - a}, \frac{d-c}{d-a} \}$.
\end{lemma}

\begin{proof}
Note that $f(t)$ is a convex function. Therefore the supremum 
of $f(t)$ is obtained in one of the interval endpoints. We see that
  \begin{align*}
  f(b) = \left|
  \frac{b-b}{c-b} - \frac{b - a}{d - a}
  \right|
  =
  \left|
  - \frac{b - a}{d - a}
  \right|
  =
  \frac{b - a}{d - a}
  \end{align*}
and
  \begin{align*}
  f(c) =
  \left|
  \frac{c-b}{c-b} - \frac{c - a}{d - a}
  \right|
  =
  \left|
  1 - \frac{c - a}{d - a}
  \right|
  =
  \left|
  1 - \frac{c - d + d- a}{d - a}
  \right|
  = \frac{d-c}{d-a}
  \end{align*}
which was what we wanted.
\end{proof}

We are now ready to show Proposition \ref{prop:innersupremum}.

\begin{proof}[Proof (of Proposition \ref{prop:innersupremum})]

We will compute the supremum over $t \in Q(\mathcal{T} \mid z) 
= [q_{\min, z}, q_{\max, z}]$ as 
the maximum of the suprema over the intervals $[q_{k, z}, q_{k+1, z}]$ for 
$k=1, \dots, m-1$, i.e., 
	\begin{align*}
	\sup_{t \in Q(\mathcal{T} \mid z)} 
  | \tilde{F}^{(m)}(t \mid z) - \hat{F}^{(m, n)}(t \mid z)| 
  & = \max_{k = 1, \dots, m-1} 
  \sup_{t \in [q_{k, z}, q_{k+1, z}]} | \tilde{F}^{(m)}(t \mid z) - \hat{F}^{(m, n)}(t \mid z) |.
	\end{align*}
This is useful since on each interval of the form 
$[q_{k, z}, q_{k+1, z}]$ we have that $\tilde{F}^{(m)}(\cdot \mid z)$ is 
a linear function, while $\hat{F}^{(m, n)}(\cdot \mid z)$ is a piecewise 
linear function.

First fix a distribution $P \in \mathcal{P}_0$ and $\epsilon > 0$.
Using Lemma 
\ref{lemma:distinctquantilepositions} we choose $N \geq 1$ such that 
$\hat{q}_{k, z} \in (q_{k-1, z}, q_{k+1, z})$ for $k=1, \dots, m-1$ and $z \in [0, 1]^d$ 
for each grid $(\tau_k)_{k=1}^m$ in $\mathcal{T}$ 
with probability at least $1 - \epsilon$. Now fix a $k=1, \dots, m-1$
such that we will examine the supremum on $[q_{k, z}, q_{k+1, z}]$. 
The relative position of the true and estimated conditional quantiles 
can be divided into four cases:
  \begin{itemize}
  \item[1)] $q_{k, z} \geq \hat{q}_{k, z}$ and $q_{k+1, z} \geq \hat{q}_{k+1, z}$. 
  \item[2)] $q_{k, z} \geq \hat{q}_{k, z}$ and $q_{k+1, z} < \hat{q}_{k+1, z}$. 
  \item[3)] $q_{k, z} < \hat{q}_{k, z}$ and $q_{k+1, z} \geq \hat{q}_{k+1, z}$. 
  \item[4)] $q_{k, z} < \hat{q}_{k, z}$ and $q_{k+1, z} < \hat{q}_{k+1, z}$.
  \end{itemize}
We start with case 1). First we compute the supremum over 
$t \in [q_{k, z}, \hat{q}_{k+1, z}]$ and then over 
$t \in [\hat{q}_{k+1, z}, q_{k+1, z}]$. We have that 
  \begin{align*}
  | \tilde{F}^{(m)} (t \mid z) - \hat{F}^{(m, n)}(t \mid z)|
  =
  (\tau_{k+1} - \tau_k) \left(
  \frac{t - \hat{q}_{k, z}}{\hat{q}_{k+1, z} - \hat{q}_{k, z}}
  -
  \frac{t - q_{k, z}}{q_{k+1, z} - q_{k, z}}
  \right)
  \end{align*}
for $t \in [q_{k, z}, \hat{q}_{k+1, z}]$. Hence we can compute the supremum as
  \begin{align*}
  & \sup_{t \in [q_{k, z}, \hat{q}_{k+1, z}]} 
  | \tilde{F}^{(m)}  (t \mid z) - \hat{F}^{(m, n)}(t \mid z)| 
  \\
  = \ &
  (\tau_{k+1} - \tau_k)
  \max \left\{
  \frac{q_{k, z} - \hat{q}_{k, z}}{\hat{q}_{k+1, z} - \hat{q}_{k, z}}
  ,
  \frac{q_{k+1, z} - \hat{q}_{k+1, z}}{q_{k+1, z} - q_{k, z}}
  \right\}
  \\
  \leq \ & 
  \kappa_m 
  \max \left\{
  \frac{\max_{k=1, \dots, m}|q_{k, z} - \hat{q}_{k, z}|}
  {\min_{k = 1, \dots, m-1}(\hat{q}_{k+1, z} - \hat{q}_{k, z})}
  ,
  \frac{\max_{k=1, \dots, m}|q_{k, z} - \hat{q}_{k, z}|}
  {\min_{k = 1, \dots, m-1}({q}_{k+1, z} - {q}_{k, z})}
  \right\}
  \end{align*}
where we have used Lemma \ref{lemma:basicinequality1}. Now we see that
  \begin{align*}
  & | \tilde{F}^{(m)} (t \mid z) - \hat{F}^{(m, n)}(t \mid z)|
  \\
  = \ &
  (\tau_{k+1} - \tau_k) + 
  (\tau_{k+2} - \tau_{k+1})
  \frac{t - \hat{q}_{k+1, z}}{\hat{q}_{k+2, z} - \hat{q}_{k+1, z}} 
  -
  (\tau_{k+1} - \tau_k) \frac{t - q_{k, z}}{q_{k+1, z} - {q}_{k, z}}
  \end{align*}
for $t \in [\hat{q}_{k+1, z}, q_{k+1, z}]$. We compute the supremum to be
  \begin{align*}
  & \sup_{t \in [\hat{q}_{k+1, z}, q_{k+1, z}]}
  | \tilde{F}^{(m)} (t \mid z) - \hat{F}^{(m, n)}(t \mid z)|
  \\
  = \ & 
  \max \left\{
  (\tau_{k+1} - \tau_k) 
  \frac{\hat{q}_{k+1, z} - q_{k+1, z}}{q_{k+1, z} - q_{k, z}},
  (\tau_{k+2} - \tau_{k+1}) 
  \frac{\hat{q}_{k+1, z} - q_{k+1, z}}{\hat{q}_{k+2, z} - \hat{q}_{k+1, z}}
  \right\}
  \\
  \leq \ & 
  \kappa_m 
  \max \left\{
  \frac{\max_{k=1, \dots, m}|q_{k, z} - \hat{q}_{k, z}|}
  {\min_{k = 1, \dots, m-1}(\hat{q}_{k+1, z} - \hat{q}_{k, z})}
  ,
  \frac{\max_{k=1, \dots, m}|q_{k, z} - \hat{q}_{k, z}|}
  {\min_{k = 1, \dots, m-1}({q}_{k+1, z} - {q}_{k, z})}
  \right\}
  \end{align*}
where we have used Lemma \ref{lemma:basicinequality2}. This covers case 1). 

Now let us proceed to case 2). Here we can evaluate the supremum over 
$t \in [q_{k, z}, q_{k+1, z}]$ directly. We have that
  \begin{align*}
  | \tilde{F}^{(m)} (t \mid z) - \hat{F}^{(m, n)}(t \mid z)|
  =
  (\tau_{k+1} - \tau_k) \left|
  \frac{t - q_{k, z}}{q_{k+1, z} - {q}_{k, z}}
  -
  \frac{t - \hat{q}_{k, z}}{\hat{q}_{k+1, z} - \hat{q}_{k, z}}
  \right|.
  \end{align*}
The supremum can now be evaluated using Lemma \ref{lemma:basicinequality3} to be
  \begin{align*}
  & \sup_{t \in [q_{k, z}, q_{k+1, z}]}
  | \tilde{F}^{(m)} (t \mid z) - \hat{F}^{(m, n)}(t \mid z)|
  \\
  = \ & 
  (\tau_{k+1} - \tau_k)
  \max \left\{
  \frac{q_{k, z} - \hat{q}_{k, z}}{\hat{q}_{k+1, z} - \hat{q}_{k, z}},
  \frac{\hat{q}_{k+1, z} - q_{k+1, z}}{\hat{q}_{k+1, z} - \hat{q}_{k, z}}
  \right\}
  \\
  \leq \ & 
  \kappa_m 
  \frac{\max_{k=1, \dots, m}|q_{k, z} - \hat{q}_{k, z}|}
  {\min_{k = 1, \dots, m-1}(\hat{q}_{k+1, z} - \hat{q}_{k, z})}. 
  \end{align*}
In case 3) we need to divide into three cases, namely when 
$t \in [q_{k, z}, \hat{q}_{z, k}]$, $t \in [\hat{q}_{k, z}, \hat{q}_{k+1, z}]$ and 
$t \in [\hat{q}_{k+1, z}, q_{k+1, z}]$. In the first case we have
  \begin{align*}
  & | \tilde{F}^{(m)} (t \mid z) - \hat{F}^{(m, n)}(t \mid z)|
  \\
  = \ &
  (\tau_k - \tau_{k-1}) + 
  (\tau_{k+1} - \tau_k)
  \frac{t - q_{k, z}}{q_{k+1, z} - {q}_{k, z}} 
  -
  (\tau_k - \tau_{k-1})
  \frac{t - \hat{q}_{k-1, z}}{\hat{q}_{k, z} - \hat{q}_{k-1, z}}
  \end{align*}
for $t \in [q_{k, z}, \hat{q}_{k, z}]$. Therefore we have
  \begin{align*}
  & \sup_{t \in [q_{k, z}, \hat{q}_{z, k}]}
  | \tilde{F}^{(m)} (t \mid z) - \hat{F}^{(m, n)}(t \mid z)|
  \\
  \leq \ & 
  \max \left\{
  (\tau_{k} - \tau_{k-1}) \frac{\hat{q}_{k, z} - q_{k, z}}{\hat{q}_{k, z} - \hat{q}_{k-1, z}},
  (\tau_{k+1} - \tau_k) \frac{\hat{q}_{k, z} - q_{k, z}}{q_{k+1, z} - {q}_{k, z}}
  \right\}
  \\
  \leq \ &
   \kappa_m 
  \max \left\{
  \frac{ \max_{k=1, \dots, m}
  |q_{k, z} - \hat{q}_{k, z}|}
  { \min_{k=1, \dots, m-1}
  (\hat{q}_{k+1, z} - \hat{q}_{k, z})}, 
  \frac{ \max_{k=1, \dots, m}
  |q_{k, z} - \hat{q}_{k, z}|}
  { \min_{k=1, \dots, m-1}
  (q_{k+1, z} - q_{k, z})}
  \right\}
  \end{align*}
where we have used Lemma \ref{lemma:basicinequality2}. In the second case we have
  \begin{align*}
  | \tilde{F}^{(m)} (t \mid z) - \hat{F}^{(m, n)}(t \mid z)|
  =
  (\tau_{k+1} - \tau_k) \left|
  \frac{t - q_{k, z}}{q_{k+1, z} - {q}_{k, z}} 
  -
  \frac{t - \hat{q}_{k, z}}{\hat{q}_{k+1, z} - \hat{q}_{k, z}}
  \right|,
  \end{align*}
for $t \in [\hat{q}_{k, z}, \hat{q}_{k+1, z}]$ and therefore we obtain 
  \begin{align*}
  & \sup_{t \in [\hat{q}_{k, z}, \hat{q}_{k+1, z}]}
  | \tilde{F}^{(m)} (t \mid z) - \hat{F}^{(m, n)}(t \mid z)|
  \\
  \leq \ & 
  (\tau_{k+1} - \tau_k)
  \max \left\{
  \frac{\hat{q}_{k, z} - q_{k, z}}{q_{k+1, z} - q_{k, z}},
  \frac{q_{k+1, z} - \hat{q}_{k+1, z}}{q_{k+1, z} - q_{k, z}}
  \right\}
  \\
  \leq \ & 
  \kappa_m \cdot \frac{ \max_{k=1, \dots, m}
  |q_{k, z} - \hat{q}_{k, z}|}
  { \min_{k=1, \dots, m-1}
  ({q}_{k+1, z} - {q}_{k, z})}
  \end{align*}
where we have used Lemma \ref{lemma:basicinequality3}. In the third case we have
  \begin{align*}
  & | \tilde{F}^{(m)} (t \mid z) - \hat{F}^{(m, n)}(t \mid z)|
  \\
  = \ & 
  (\tau_{k+1} - \tau_k) + 
  (\tau_{k+2} - \tau_{k+1})
  \frac{t - \hat{q}_{k+1, z}}{\hat{q}_{k+2, z} - \hat{q}_{k+1, z}} 
  -
  (\tau_{k+1} - \tau_k) \frac{t - q_{k, z}}{q_{k+1, z} - {q}_{k, z}}
  \end{align*}
for $t \in [\hat{q}_{k+1, z}, q_{k+1, z}]$. So we obtain 
  \begin{align*}
  & \sup_{t \in [\hat{q}_{k+1, z}, q_{k+1, z}]}
  | \tilde{F}^{(m)} (t \mid z) - \hat{F}^{(m, n)}(t \mid z)|
  \\
  = \ & 
  \max \left\{
  (\tau_{k+1} - \tau_k) 
  \frac{\hat{q}_{k+1, z} - q_{k+1, z}}{q_{k+1, z} - q_{k, z}},
  (\tau_{k+2} - \tau_{k+1}) 
  \frac{\hat{q}_{k+1, z} - q_{k+1, z}}{\hat{q}_{k+2, z} - \hat{q}_{k+1, z}}
  \right\}
  \\
  \leq \ & 
  \kappa_m 
  \max \left\{
  \frac{\max_{k=1, \dots, m}|q_{k, z} - \hat{q}_{k, z}|}
  {\min_{k = 1, \dots, m-1}(\hat{q}_{k+1, z} - \hat{q}_{k, z})}
  ,
  \frac{\max_{k=1, \dots, m}|q_{k, z} - \hat{q}_{k, z}|}
  {\min_{k = 1, \dots, m-1}({q}_{k+1, z} - {q}_{k, z})}
  \right\}
  \end{align*}
where we have used Lemma \ref{lemma:basicinequality2}.

Let us now examine case 4). Here we have the two sub cases 
$t \in [q_{k, z}, \hat{q}_{k, z}]$ and $t \in [\hat{q}_{k, z}, q_{k+1, z}]$. 
First we see that 
  \begin{align*}
  & | \tilde{F}^{(m)} (t \mid z) - \hat{F}^{(m, n)}(t \mid z)|
  \\
  = & \
  (\tau_k - \tau_{k-1}) + 
  (\tau_{k+1} - \tau_k)
  \frac{t - q_{k, z}}{q_{k+1, z} - {q}_{k, z}} 
  -
  (\tau_k - \tau_{k-1})
  \frac{t - \hat{q}_{k-1, z}}{\hat{q}_{k, z} - \hat{q}_{k-1, z}}
  \end{align*}
for $t \in [q_{k, z}, \hat{q}_{k, z}]$. Thus we have 
  \begin{align*}
  & \sup_{t \in [q_{k, z}, \hat{q}_{z, k}]}
  | \tilde{F}^{(m)} (t \mid z) - \hat{F}^{(m, n)}(t \mid z)|
  \\
  \leq \ & 
  \max \left\{
  (\tau_{k} - \tau_{k-1}) \frac{\hat{q}_{k, z} - q_{k, z}}{\hat{q}_{k, z} - \hat{q}_{k-1, z}},
  (\tau_{k+1} - \tau_k) \frac{\hat{q}_{k, z} - q_{k, z}}{q_{k+1, z} - {q}_{k, z}}
  \right\}
  \\
  \leq \ &
   \kappa_m 
  \max \left\{
  \frac{ \max_{k=1, \dots, m}
  |q_{k, z} - \hat{q}_{k, z}|}
  { \min_{k=1, \dots, m-1}
  (\hat{q}_{k+1, z} - \hat{q}_{k, z})}, 
  \frac{ \max_{k=1, \dots, m}
  |q_{k, z} - \hat{q}_{k, z}|}
  { \min_{k=1, \dots, m-1}
  (q_{k+1, z} - q_{k, z})}
  \right\}
  \end{align*}
where we have used \ref{lemma:basicinequality2}. Now in the second case we have
  \begin{align*}
  | \tilde{F}^{(m)} (t \mid z) - \hat{F}^{(m, n)}(t \mid z)|
  =
  (\tau_{k+1} - \tau_k) \left(
  \frac{t - q_{k, z}}{q_{k+1, z} - {q}_{k, z}} 
  -
  \frac{t - \hat{q}_{k, z}}{\hat{q}_{k+1, z} - \hat{q}_{k, z}}
  \right)
  \end{align*}
for $t \in [\hat{q}_{k, z}, q_{k+1, z}]$. From this we get the supremum to be
  \begin{align*}
  & \sup_{t \in [\hat{q}_{k, z}, q_{k+1, z}]}
  | \tilde{F}^{(m)} (t \mid z) - \hat{F}^{(m, n)}(t \mid z)|
  \\
  = \ & 
  (\tau_{k+1} - \tau_k)
  \max
  \left\{
  \frac{\hat{q}_{k, z} - q_{k, z}}{q_{k+1, z} - q_{k, z}},
  \frac{q_{k+1, z} - \hat{q}_{k+1, z}}{\hat{q}_{k+1, z}-\hat{q}_{k, z}}
  \right\}
  \\
  \leq \ &
   \kappa_m 
  \max \left\{
  \frac{ \max_{k=1, \dots, m}
  |q_{k, z} - \hat{q}_{k, z}|}
  { \min_{k=1, \dots, m-1}
  (\hat{q}_{k+1, z} - \hat{q}_{k, z})}, 
  \frac{ \max_{k=1, \dots, m}
  |q_{k, z} - \hat{q}_{k, z}|}
  { \min_{k=1, \dots, m-1}
  (q_{k+1, z} - q_{k, z})}
  \right\}
  \end{align*}
where we have used Lemma \ref{lemma:basicinequality1}. 
Taking maximum of all cases and sub cases yields the desired result.
\end{proof}

We will now move on to tackling the problem of controlling the outer 
supremum over $z \in [0, 1]^d$. First we prove the following technical lemma 
that gives control over the denominators in $A(m, n, z)$ and $B(m, n, z)$. 

\begin{lemma} \label{lemma:controllingdenominator}
Let Assumption \ref{assump:Qhatconsistency} be satisfied. Let 
$\gamma_m := \min_{k = 1, \dots, m-1}(\tau_{k+1} - \tau_k)$ denote the finest subinterval of the grid. 
Then for each $P \in \mathcal{P}_0$ we have
  \begin{align*}
  \min_{k=1, \dots, m-1}(q_{k+1, z} - q_{k, z})
  \geq \frac{\gamma_m}{C_P}
  \end{align*}
for almost all $z \in [0, 1]^d$ for each grid 
$(\tau_k)_{k=1}^m$ in $\mathcal{T}$. Also for 
all $\epsilon > 0$ there is $N \geq 1$ such that for all $n \geq N$ we
have
  \begin{align*}
  \min_{k=1, \dots, m-1}(\hat{q}_{k+1, z} - \hat{q}_{k, z})
  \geq \frac{\gamma_m}{3C_P}
  \end{align*}
for almost all $z \in [0, 1]^d$ for each grid $(\tau_k)_{k=1}^m$ in 
$\mathcal{T}$ with probability at least 
$1 - \epsilon$.
\end{lemma}

\begin{proof}
Fix a distribution $P \in \mathcal{P}_0$. We see that
  \begin{align*}
  \tau_{k+1} - \tau_k 
  & =
  F(q_{k+1, z} \mid z) - F(q_{k, z} \mid z)
  \\
  & =
  \int_{q_{k, z}}^{q_{k+1, z}} f(x \mid z) \mathrm{d} x
  \leq C_P \cdot (q_{k+1, z} - q_{k, z})
  \end{align*}
for each $k=1, \dots, m-1$ and almost all $z \in [0, 1]^d$ for each grid 
$(\tau_k)_{k=1}^m$ in $\mathcal{T}$. Here we have 
used Assumption \ref{assump:Qhatconsistency} (ii).
Rearranging and taking minimum, we have that 
  \begin{align*}
  \min_{k=1, \dots, m-1} (q_{k+1, z} - q_{k, z})
  \geq
  \min_{k=1, \dots, m-1}
  \frac{\tau_{k+1} - \tau_k}{C_P} = \frac{\gamma_m}{C_P}
  \end{align*}
for almost all $z \in [0, 1]^d$ and each grid 
$(\tau_k)_{k=1}^m$ in $\mathcal{T}$. 
Now let $\epsilon > 0$ be given. Choose $N \geq 1$ such that 
for all $n \geq N$ we have 
  \begin{align*}
  P \left( \hat{q}_{k, z} \in 
  \left(
  q_{k, z} - \frac{\gamma_m}{3 C_P}, q_{k, z} + \frac{\gamma_m}{3 C}
  \right) \right) \geq 1 - \epsilon
  \end{align*}
for all $k=1, \dots, m$ and all $z \in [0, 1]^d$ for each 
$(\tau_k)_{k=1}^m$ in $\mathcal{T}$, which is possible due to 
Assumption \ref{assump:Qhatconsistency} (i). In this case
  \begin{align*}
  \hat{q}_{k, z} \leq q_{k, z} + \frac{\gamma_m}{3 C_P}
  \quad
  \text{and}
  \quad
  \hat{q}_{k+1, z} \geq q_{k+1, z} -  \frac{\gamma_m}{3 C_P}
  \end{align*}
for all $k=1, \dots, m-1$ and $z \in [0, 1]^d$ with probability 
at least $1 - \epsilon$. Thus for $n \geq N$,
  \begin{align*}
  \min_{k=1, \dots, m-1} (\hat{q}_{k+1, z} - \hat{q}_{k, z})
  & \geq 
  \min_{k=1, \dots, m-1} \left(
  q_{k+1, z} -  \frac{\gamma_m}{3 C_P}
  -
  \left(
  q_{k, z} +  \frac{\gamma_m}{3 C_P}
  \right)
  \right)
  \\
  & = 
  \min_{k=1, \dots, m-1} 
  \left(
  q_{k+1, z} - q_{k, z} -  \frac{2 \gamma_m}{3 C_P}
  \right)
  \\
  & \geq  \frac{\gamma_m}{C_P} -  \frac{2\gamma_m}{3 C_P}
  =
   \frac{\gamma_m}{3C_P}
  \end{align*}
for all $z \in [0, 1]^d$ and each grid 
$(\tau_k)_{k=1}^m$ in $\mathcal{T}$ with 
probability at least $1 - \epsilon$.
\end{proof}

We are now ready to prove the main result.

\begin{proof}[Proof (of Proposition \ref{prop:Fhatconsistent})]

Fix a distribution $P \in \mathcal{P}_0$. Let $\epsilon \in (0, 1)$ be given. 
Firstly, we use Proposition \ref{prop:innersupremum} to choose $N_1 \geq 1$ such that the 
event 
  \begin{align*}
  E_1 =
  \left(
  \sup_{t \in Q(\mathcal{T} \mid z)} |\tilde{F}^{(m)}(t \mid z) - \hat{F}^{(m, n)}(t \mid z)|
  \leq \max \{ A(m, n, z), B(m, n, z) \}
  \right)
  \end{align*}
has probability at least $1 - \epsilon / 3$ for all $n \geq N_1$ and every grid 
$(\tau_k)_{k=1}^m$ in $\mathcal{T}$. Secondly, 
according to Lemma \ref{lemma:controllingdenominator} we have that
  \begin{align*}
  \min_{k=1, \dots, m-1}(q_{k+1, z} - q_{k, z})
  \geq \frac{\gamma_m}{C_P}
  \end{align*}
and we can choose $N_2 \geq 1$ such that the event 
  \begin{align*}
  E_2 =  \left(
  \min_{k=1, \dots, m-1}(\hat{q}_{k+1, z} - \hat{q}_{k, z}) \geq 
  \frac{\gamma_m}{3C_P}
  \right)
  \end{align*}
has probability at least $1 - \epsilon / 3$ for all $n \geq N_2$ and every grid 
$(\tau_k)_{k=1}^m$ in $\mathcal{T}$. Thirdly, 
we can choose $N_3 \geq 1$ and $M_P' > 0$ such that the event 
  \begin{align*}
  E_3 = \left( 
  \frac{\sup_{z \in [0, 1]^d} \max_{k=1, \dots, m} |q_{k, z} - \hat{q}_{k, z}|}
  {g_P(n)} \leq M'_P
  \right)
  \end{align*}
has probability at least $1 - \epsilon / 3$ for all $n \geq N_3$ and every
$(\tau_k)_{k=1}^m$ in $\mathcal{T}$ using Assumption \ref{assump:Qhatconsistency} (i). 
Now we note that on the event $E := E_1 \cap E_2 \cap E_3$ we have
  \begin{align*}
  \|\tilde{F}^{(m)} - \hat{F}^{(m, n)}\|_{\mathcal{T}, \infty}
  & \leq 
  \sup_{z \in [0, 1]^d} \max \{ A(m, n, z), B(m, n, z) \}
  \\
  & = 3 C_P \cdot \frac{\kappa_m}{\gamma_m} 
  \sup_{z \in [0, 1]^d} \max_{k=1, \dots, m} |q_{k, z} - \hat{q}_{k, z}|
  \\
  & \leq 3 C_P \cdot M_P' \cdot g_P(n)
  \end{align*}
with probability $P(E) \geq 1 - \epsilon$ for all $n \geq N$ and every grid 
$(\tau_k)_{k=1}^m$ in $\mathcal{T}$ where $N := \max\{N_1, N_2, N_3\}$. Here we have 
used that $\kappa_m / \gamma_m = 1$ due to the grids being equidistant.
We can now set $M_P := 3 C_P \cdot M_P'$ such that
  \begin{align*}
  P \left( 
  \frac{\|\tilde{F}^{(m)} - \hat{F}^{(m, n)}\|_{\mathcal{T}, \infty}}{g_P(n)} > M_P
  \right) < \epsilon
  \end{align*}
whenever $n \geq N$. This shows that 
$\|\tilde{F}^{(m)} - \hat{F}^{(m, n)}\|_{\mathcal{T}, \infty} \in \mathcal{O}_P(g_P(n))$ for every 
equidistant grid $(\tau_k)_{k=1}^m$ in $\mathcal{T}$ as wanted. 
\end{proof}

\subsection{Proof of Theorem \ref{thm:pointwiseconsistency}}

According to Corollary \ref{cor:keyinequality} we have
  \begin{align*}
  \|F -  \hat{F}^{(m_n, n)}\|_\infty
  \leq 
  \kappa_{m_n}
  +
  \|\tilde{F}^{(m_n)} - \hat{F}^{(m_n, n)}\|_\infty.
  \end{align*}
Here $\|\tilde{F}^{(m_n)} - \hat{F}^{(m_n, n)}\|_\infty \in \mathcal{O}_P(g_P(n))$ for 
each equidistant grid $(\tau_k)_{k=1}^{m_n}$ in $\mathcal{T}$ due to Proposition \ref{prop:Fhatconsistent}. 
Since we have assumed that $\kappa_{m_n} \in o(g_P(n))$ we have the result. 
\hfill $\square$

\subsection{Proof of Proposition \ref{prop:uniformFhatconsistent}}

The proof follows immediately from the proof of Proposition \ref{prop:Fhatconsistent} and 
the stronger Assumption \ref{assump:uniformQhatconsistent} in the following way. 
Note that the statement of Lemma \ref{lemma:distinctquantilepositions} holds uniformly over 
$P \in \mathcal{P}_0$ under Assumption \ref{assump:uniformQhatconsistent} (i). Therefore 
Proposition \ref{prop:innersupremum} also holds uniformly over $P \in \mathcal{P}_0$. Furthermore, 
the result of Lemma \ref{lemma:controllingdenominator} also holds uniformly in $P \in \mathcal{P}_0$ 
under Assumption \ref{assump:uniformQhatconsistent}. Therefore the probability of the events $E_1, E_2$ and 
$E_3$ can be controlled uniformly over $P \in \mathcal{P}_0$ from which the result follows.
\hfill $\square$

\subsection{Proof of Theorem \ref{thm:uniformconsistency}}

The corollary follows from Proposition \ref{prop:uniformFhatconsistent} using the same 
argument as in the proof of Theorem \ref{thm:pointwiseconsistency}. \hfill $\square$

\subsection{Proof of Corollary \ref{cor:Qhatconsistency}}

Using Theorem \ref{thm:belloni} we have that 
  \begin{align*}
  \sup_{z \in [0, 1]^d} \sup_{\tau \in \mathcal{Q}} |Q(\tau \mid z) - \hat{Q}(\tau \mid z)| 
  & = 
  \sup_{z \in [0, 1]^d} \sup_{\tau \in \mathcal{Q}} |h(z)^T (\beta_\tau - \hat{\beta}_\tau)|
  \\
  & \leq \sup_{z \in [0, 1]^d} \|h(z)\|_2 
  \sup_{\tau \in \mathcal{Q}} \|\beta_\tau - \hat{\beta}_\tau\|_2 
  \\
  & \in \mathcal{O}_P \left( 
  \sqrt{\frac{s_n \log(p \vee n)}{n}}
  \right)
  \end{align*}
since $\sup_{z \in [0, 1]^d} \|h(z)\|_2 < \infty$ because $[0,1]^d$ is compact and $h$ is continuous. \hfill $\square$

\subsection{Proof of Proposition \ref{prop:XYZtoUU}}
Assume that $X \indep Y \mid Z$. Then it also holds that 
$(X, Z) \indep (Y, Z) \indep Z$ and thus $U_1 \indep U_2 \mid Z$. 
Letting $f$ denote a generic density function, we now have that
  \begin{align*}
  f(u_1, u_2) & = \int f(u_1, u_2 \mid z) f(z) \mathrm{d} z
  =
  \int f(u_1 \mid z) f(u_2 \mid z) f(z) \mathrm{d} z
  \\
  & = \int f(u_1) f(u_2) f(z) \mathrm{d} z 
  =
  f(u_1) f(u_2)
  \end{align*}
for all $u_1, u_2 \in [0, 1]$, where we have used Proposition 
\ref{prop:Uuniform}. \hfill $\square$

\subsection{Proof of Theorem \ref{thm:rhoasympnormal}}

Before proving the theorem, we will supply a lemma  
that will aid us during the proof. 

\begin{lemma} \label{lem:mylittlelemma}
Let $\hat{F}_{X \mid Z}^{(n)}$ and $\hat{F}^{(n)}_{Y \mid Z}$ 
satisfy Assumption \ref{assump:Fhatconsistency}. Then 
	\begin{align*}
	\| \phi_k \circ F_{X \mid Z} - \phi_k \circ \hat{F}_{X \mid Z}^{(n)} \|_\infty
	\in \mathcal{O}_P(g_P(n))
	\quad \text{and} \quad
	\| \phi_k \circ F_{Y \mid Z} - \phi_k \circ \hat{F}_{Y \mid Z}^{(n)} \|_\infty 
	\in \mathcal{O}_P(h_P(n))
	\end{align*}
for each $k=1, \dots, q$ given that $\phi$ satisfies 
Assumption \ref{assump:phimap}.
\end{lemma}

\begin{proof}
We only show the first statement. Fix $k=1, \dots, q$. 
We need to control the supremum
	\begin{align*}
	\sup_{z \in [0, 1]^d} \sup_{t \in [0, 1]} 
	|
	\phi_k (F_{X \mid Z}(t \mid z)) - 
	\phi_k (\hat{F}_{X \mid Z}^{(n)}(t \mid z))
	|.
	\end{align*}
We will divide the supremum over $t \in [0, 1]$ into two cases. Namely, 
when $t \in Q(\mathcal{T} \mid z) = [q_{\min, z}, q_{\max, z}]$ and when $t \in 
Q(\mathcal{T}^c \mid z) = [q_{\min, z}, q_{\max, z}]^c$. First we see that
	\begin{align*}
	\sup_{z \in [0, 1]^d} \sup_{t \in Q(\mathcal{T} \mid z)} 
	|
	\phi_k & (F_{X \mid Z}(t \mid z)) - 
	\phi_k (\hat{F}_{X \mid Z}^{(n)}(t \mid z))
	|
	\\
	& \leq L_k \cdot \|F_{X \mid Z} - 
	\hat{F}^{(n)}_{X \mid Z}\|_{\mathcal{T}, \infty} 
	\in \mathcal{O}_P(g_P(n))
	\end{align*}
where $L_k$ is the Lipschitz constant of $\phi_k$ under Assumption 
\ref{assump:phimap} (ii). Here we have used the consistency in
Assumption \ref{assump:Fhatconsistency} (i). Next we examine the 
supremum over $t \in Q(\mathcal{T}^c \mid z)$. First note that 
$F_{X \mid Z}(t \mid z) \in [\tau_{\min}, \tau_{\max}]^c$ whenever 
$t \in Q(\mathcal{T}^c \mid z)$. Also recall that the support 
of $\phi_k$ is 
$\mathcal{T}_k \subset \mathcal{T} = [\tau_{\min}, \tau_{\max}]$. 
Therefore $\phi_k(F_{X \mid Z}(t \mid z)) = 0$ for 
$t \in Q(\mathcal{T}^c \mid z)$. Hence we have
	\begin{align*}
	\sup_{z \in [0, 1]^d} \sup_{t \in Q(\mathcal{T}^c \mid z)} 
	|
	\phi_k ( F_{X \mid Z} (t \mid z) ) -
	\phi_k ( \hat{F}_{X \mid Z}^{(n)} (t \mid z) )
	| = 
	\sup_{z \in [0, 1]^d} \sup_{t \in Q(\mathcal{T}^c \mid z)} 
	|
	\phi_k ( \hat{F}_{X \mid Z}^{(n)}(t \mid z) )
	|.
	\end{align*}
By Assumption \ref{assump:Fhatconsistency} (i) we know that 
	\begin{align*}
	\hat{F}^{(n)}_{X \mid Z}(q_{\min, z} \mid z) \stackrel{P}{\to} \tau_{\min}
	\quad \text{and} \quad 
	\hat{F}^{(n)}_{X \mid Z}(q_{\max, z} \mid z) \stackrel{P}{\to} \tau_{\max}
	\end{align*}
for all $z \in [0, 1]^d$. Since $\hat{F}^{(n)}_{X \mid Z}(\cdot \mid z)$ is increasing 
we thus know that the limit $\xi(t, z)$ from Assumption \ref{assump:Fhatconsistency} (ii) 
must satisfy 
$\xi(t, z) \in [\tau_{\min}, \tau_{\max}]^c$ for $t \in Q(\mathcal{T}^c \mid z)$ and 
$z \in [0, 1]^d$. Again, since the support of $\phi_k$ is 
$\mathcal{T}_k \subset \mathcal{T} = [\tau_{\min}, \tau_{\max}]$ we have that 
$\phi_k(\xi(t, z)) = 0$ when $t \in Q(\mathcal{T}^c \mid z)$ and $z \in [0, 1]^d$. 
Therefore we have that
	\begin{align*}
	\sup_{z \in [0, 1]^d} \sup_{t \in Q(\mathcal{T}^c \mid z)} 
	|
	\phi_k ( \hat{F}_{X \mid Z}^{(n)}(t \mid z) )
	| & = 
	\sup_{z \in [0, 1]^d} \sup_{t \in Q(\mathcal{T}^c \mid z)} 
	|\phi_k(\xi(t, z))
	-
	\phi_k ( \hat{F}_{X \mid Z}^{(n)}(t \mid z) )
	| 
	\\
	& \leq L_k \cdot \| \xi - \hat{F}^{(n)}_{X \mid Z} \|_{\mathcal{T}^c, \infty} 
	\in \mathcal{O}_P(g_P(n)), 
	\end{align*}
where we have used Assumption \ref{assump:Fhatconsistency} (ii). 
Putting the two cases together we have that
	\begin{align*}
	\| \phi_k \circ F_{X \mid Z} - \phi_k \circ \hat{F}_{X \mid Z}^{(n)} \|_\infty
	\in \mathcal{O}_P(g_P(n))
	\end{align*}
which was what we wanted.
\end{proof}

We can now prove the main theorem.

\begin{proof}[Proof (of Theorem \ref{thm:rhoasympnormal})]
Fix a distribution $P \in \mathcal{H}_0$.
The key to proving the theorem is the decomposition
	\begin{align*}
	\hat{\rho}_n = \alpha_n + \beta_n + \gamma_n + \delta_n
	\end{align*}
where $\alpha_n, \beta_n, \gamma_n$ 
and $\delta_n$ are given by
	\begin{align*}
	\alpha_n & = 
    \frac{1}{n} \sum_{i=1}^n \phi(U_{1, i}) \phi(U_{2, i})^T, 
    \\
    \beta_n & = \frac{1}{n} \sum_{i=1}^n
    \left(
    \phi(\hat{U}_{1, i}) - \phi(U_{1, i})
    \right)
    \left(
    \phi(\hat{U}_{2, i}) - \phi(U_{2, i})
    \right)^T
    \\
    \gamma_n & = \frac{1}{n} \sum_{i=1}^n 
    \phi(U_{1, i}) \left( 
    \phi(\hat{U}_{2, i}) - \phi(U_{2, i})
    \right)^T, 
    \\
    \delta_n & = \frac{1}{n} \sum_{i=1}^n 
    \left( 
    \phi(\hat{U}_{1, i}) - \phi(U_{1, i})
    \right) \phi(U_{2, i})^T
	\end{align*}
The term $\alpha_n$ will be driving the asymptotics of 
the test statistics, while $\beta_n, \gamma_n$ and $\delta_n$ are 
error terms that we wish to show converge to zero sufficiently fast.

Let us start by examining $\alpha_n$. 
Under Assumption \ref{assump:phimap} (iii) we see that 
	\begin{align*}
	E_P(\phi(U_{1, i}) \phi(U_{2, i})^T) = 
	E_P(\phi(U_{1, i})) E_P( \phi(U_{2, i}))^T
	= 0
	\end{align*}
because $P \in \mathcal{H}_0$ and furthermore we see that 
	\begin{align*}
	\text{Cov}_P(\phi_k(U_{1, i}) \phi_\ell(U_{2, i}), 
	\phi_s(U_{1, i}) \phi_t(U_{2, i}))
	& =E_P(\phi_k(U_{1, i}) \phi_\ell(U_{2, i})
	\phi_s(U_{1, i}) \phi_t(U_{2, i}))
	\\
	& =
	E_P(\phi_k(U_{1, i})\phi_s(U_{1, i})) 
	E_P(\phi_\ell(U_{2, i})\phi_t(U_{2, i}))
	\\
	& =
	\int_0^1 \phi_k(u) \phi_s(u) \mathrm{d}u 
	\int_0^1 \phi_\ell(u) \phi_t(u) \mathrm{d}u 
	\\
	& = 
	\Sigma_{ks} \Sigma_{\ell t} = (\Sigma \otimes \Sigma)_{k \ell, st}
	\end{align*}
for $k, \ell, s, t = 1, \dots, q$. Observe that 
$\Sigma_{k, k} = 1$. Since $\alpha_n$ is the average 
of i.i.d. terms with zero mean and covariance $\Sigma \otimes \Sigma$, the central 
limit theorem states that
	\begin{align*}
	\sqrt{n} \alpha_n \Rightarrow_P \mathcal{N}(0, \Sigma \otimes \Sigma)
	\end{align*}
for each $P \in \mathcal{H}_0$. 

Now let us examine the term 
$\sqrt{n} \beta_n$. Fix $k, \ell = 1, \dots, q$. Then we have
	\begin{align*}
	 | \sqrt{n} \beta_{k \ell, n}| & \leq
	\frac{1}{\sqrt{n}} \sum_{i = 1}^n 
	\left|
	\phi_k (\hat{U}_{1, i}) - \phi_k (U_{1, i})
	\right| \cdot 
	\left|
	\phi_\ell (\hat{U}_{2, i}) - \phi_\ell (U_{2, i})
	\right|
	\\
	& \leq \frac{n}{\sqrt{n}}
	\|
	\phi_k \circ \hat{F}^{(n)}_{X \mid Z} 
	- 
	\phi_k \circ F_{X \mid Z}
	\|_{\infty} \cdot
	\|
	\phi_\ell \circ \hat{F}^{(n)}_{Y \mid Z} 
	-
	\phi_\ell \circ F_{Y \mid Z}
	\|_{\infty} 
	\\
	& \in \mathcal{O}_P(\sqrt{n} g_P(n) h_P(n))
	\end{align*}
where we have used Lemma \ref{lem:mylittlelemma}, which is valid due to 
Assumption \ref{assump:Fhatconsistency}. Since we have assumed that the 
rate functions satisfy $\sqrt{n} g_P(n) h_P(n) \to 0$ we can conclude that 
$| \sqrt{n} \beta_{k \ell, n}	| \to_P 0$ for each $k, \ell = 1, \dots, q$. 
Hence $\sqrt{n} \beta_n  \to_P 0$. 

Now we turn to the cross terms $\gamma_n$ and $\delta_n$. The two terms 
are dealt with analogously, so we only examine $\gamma_n$. Fix 
$k, \ell = 1, \dots, q$ and consider writing
	\begin{align*}
	\gamma_{k \ell, n} = \frac{1}{n} \sum_{i = 1}^n C_i 
	\quad \text{where} \quad 
	C_i = \phi_k(U_{1, i}) \left(
	\phi_\ell ( \hat{U}_{2, i} ) - \phi_\ell (U_{2, i})
	\right).
	\end{align*}
We will compute the mean and variance of $\sqrt{n} \gamma_{k \ell, n}$ 
conditionally on $(Y_j, Z_j)_{j = 1}^n$ in order to use Chebyshev's 
inequality to show that it converges to zero in probability. 
Observe that
	\begin{align*}
	E_P(C_i \mid (Y_j, Z_j)_{j = 1}^n) & = 
  E_P \left(
  \phi_k (U_{1, i}) 
  \left( 
  \phi_\ell (\hat{U}_{2, i}) - \phi_\ell (U_{2, i})
  \right) \mid (Y_j, Z_j)_{j = 1}^n
  \right)
  \\
  & = 
  \left( 
  \phi_\ell (\hat{U}_{2, i}) - \phi_\ell (U_{2, i})
  \right)
  E_P \left( \phi_k (U_{1, i}) \mid (Y_j, Z_j)_{j = 1}^n \right)
   \quad \mathrm{a.s.}
	\end{align*}
Here we have exploited that $\phi_\ell (U_{2, i})$ and 
$\phi_\ell(\hat{U}_{2, i}) = \phi_\ell(\hat{F}_{Y \mid Z}(Y_i \mid Z_i))$ are measurable functions of $(Y_j, Z_j)_{j = 1}^n$. 
Now since $P \in \mathcal{H}_0$ we have $\phi_k(U_{1, i}) \indep Y_i \mid Z_i$ and 
$\phi_k(U_{1, i}) \indep Z_i$ due to Proposition \ref{prop:Uuniform}. Therefore  
  \begin{align*}
  E_P \left( \phi_k (U_{1, i}) \mid (Y_j, Z_j)_{j = 1}^n \right)
  = E_P \left( 
  \phi_k (U_{1, i})
  \right) = 0 \quad \mathrm{a.s.}
  \end{align*}
where we have used Assumption \ref{assump:phimap} (iii). Hence $E_P(C_i \mid (Y_j, Z_j)_{j = 1}^n) = 0$ a.s.
From the tower property we also obtain that $E_P(C_i)=0$
and therefore $\sqrt{n} \gamma_{k \ell, n}$ has mean zero.
Let us turn to the conditional variance.
Conditionally on $(Y_j, Z_j)_{j = 1}^n$ the terms $(C_i)_{i=1}^n$ are i.i.d. because 
$\phi_\ell \circ \hat{F}_{Y \mid Z}$ is $(Y_j, Z_j)_{j = 1}^n$-measurable 
as exploited before. So we have
  \begin{align*}
  V_P(\sqrt{n} \gamma_{k \ell, n} \mid (Y_j, Z_j)_{j = 1}^n ) 
  & = \frac{1}{n} \sum_{i = 1}^n
  V_P( C_i \mid (Y_j, Z_j)_{j = 1}^n)
  =
  V_P( C_i \mid (Y_j, Z_j)_{j = 1}^n).
  \end{align*}
We compute the conditional variance to be
  \begin{align*}
  V_P( C_i \mid (Y_j, Z_j)_{j = 1}^n)
  & =
  E_P \left(
  \phi_k (U_{1, i})^2
  \left( 
  \phi_\ell (\hat{U}_{2, i}) - \phi_\ell (U_{2, i})
  \right)^2
  \mid (Y_j, Z_j)_{j = 1}^n
  \right)
  \\
  & = 
  \left( 
  \phi_\ell (\hat{U}_{2, i}) - \phi_\ell (U_{2, i})
  \right)^2 E_P \left(
  \phi_k (U_{1, i})^2 \mid (Y_j, Z_j)_{j = 1}^n
  \right)
  \\
  & = 
  \left( 
  \phi_\ell (\hat{U}_{2, i}) - \phi_\ell (U_{2, i})
  \right)^2 E_P \left(
  \phi_k (U_{1, i})^2
  \right) \\
  & = 
  \left( 
  \phi_\ell (\hat{U}_{2, i}) - \phi_\ell (U_{2, i})
  \right)^2 \quad \mathrm{a.s.}
  \end{align*}
where we have used Assumption \ref{assump:phimap} (iii). 
We can use the the law of 
total variance to see that
  \begin{align*}
  V_P(\sqrt{n} \gamma_{k \ell, n})
  & =
  E_P (V_P(\sqrt{n} \gamma_{k \ell, n} \mid (Y_j, Z_j)_{j = 1}^n))
  +
  V_P(E_P(\sqrt{n} \gamma_{k \ell, n} \mid (Y_j, Z_j)_{j = 1}^n))
  \\
  & = E_P \left( 
  \phi_\ell (\hat{U}_{2, i}) - \phi_\ell (U_{2, i})
  \right)^2 + 0 = E_P \left( 
  \phi_\ell (\hat{U}_{2, i}) - \phi_\ell (U_{2, i})
  \right)^2. 
  \end{align*}
By Lemma \ref{lem:mylittlelemma} we have that 
$\left(\phi_\ell (\hat{U}_{2, i}) - \phi_\ell (U_{2, i})\right)^2 \to_P 0$ with similar arguments as
before. Note that $\phi_\ell : [0, 1] \to \R$ is bounded due to continuity of $\phi_\ell$ and compactness of $[0, 1]$. 
Hence each term in the sequence 
  \begin{align*}
  \left(\left(\phi_\ell (\hat{U}_{2, i}) - \phi_\ell (U_{2, i}
  )\right)^2\right)_{i=1, \dots, n}
  \end{align*}
is bounded. Therefore we also have $E_P \left(\phi_\ell (\hat{U}_{2, i}) - \phi_\ell (U_{2, i})\right)^2
\to 0$. For given $\epsilon > 0$ we have by Chebyshev's inequality that
	\begin{align*}
	P(|\sqrt{n} \gamma_{k \ell, n}| > \epsilon) \leq 
	  \frac{V_P(\sqrt{n} \gamma_{k \ell, n})}{\epsilon^2}
	  =
	  \frac{1}{\epsilon^2} \cdot 
	  E_P \left(\phi_\ell (\hat{U}_{2, i}) - \phi_\ell (U_{2, i})\right)^2 \to 0
	\end{align*}
for each $P \in \mathcal{H}_0$. This shows 
$\sqrt{n} \gamma_n \to_P 0$. By the same
argument it can be shown that \mbox{$\sqrt{n} \delta_n \to_P 0$}. 
By Slutsky's lemma we now have that
  \begin{align*}
  \sqrt{n} \hat{\rho}_n = 
  \sqrt{n} \alpha_n + \sqrt{n} \beta_n 
  + \sqrt{n} \gamma_n + \sqrt{n} \delta_n 
  \ \Rightarrow_P \ \mathcal{N}(0, \Sigma \otimes \Sigma )
  \end{align*}
for each $P \in \mathcal{H}_0$. This shows the theorem.
\end{proof}

\subsection{Proof of Corollary \ref{cor:chisqteststat}}

First note that $\Sigma$ is a positive definite matrix as $\phi_1, \ldots, \phi_q$
are assumed linearly independent. It thus has a positive definite  
matrix square root $\Sigma^{- 1 / 2}$ satisfying 
$\Sigma^{-1/2} \Sigma \Sigma^{-1/2} = I$, and we have that
	\begin{align*}
	\sqrt{n} \Sigma^{-1 / 2} \hat{\rho}_n \Sigma^{-1 / 2} \Rightarrow_P
	\mathcal{N}(0, I \otimes I)
	\end{align*}
for $P \in \mathcal{H}_0$ where we have used 
Theorem \ref{thm:rhohat_to_rho}. The test statistics $T_n$ is therefore 
well defined and
	\begin{align*}
	n T_n  = \| \sqrt{n} \Sigma^{-1 / 2} \hat{\rho}_n \Sigma^{-1 / 2} \|_F^2
	 \Rightarrow_P \chi^2_{q^2}
	\end{align*}
for $P \in \mathcal{H}_0$ by the continuous mapping theorem. \hfill $\square$

\subsection{Proof of Corollary \ref{cor:level}}

Under Assumption \ref{assump:Fhatconsistency} we have by 
Corollary \ref{cor:chisqteststat} that 
$n T_n \Rightarrow_P \chi^2_{q^2}$. Therefore
	\begin{align*}
	\limsup_{n \to \infty} 
    E_P(\hat{\Psi}_n) & = \limsup_{n \to \infty} 
    P( n T_n > z_{1 - \alpha}) 
    = \limsup_{n \to \infty} 
    (1 - (F_{n T_n}(z_{1 - \alpha})) = 1 - (1 - \alpha) = \alpha. 
	\end{align*}
because $F_{n T_n}(t) \to \Phi(t)$ as $n \to \infty$ for 
all $t \in \R$ where $\Phi$ is the distribution function 
of a $\chi^2_{q^2}$-distribution and $z_{1 - \alpha}$ is the 
$(1 - \alpha)$-quantile of a $\chi^2_{q^2}$-distribution. \hfill $\square$

\subsection{Proof of Theorem \ref{thm:rhohat_to_rho}}

The proof uses the same decomposition as in the proof of Theorem \ref{thm:rhoasympnormal},
i.e., $\hat{\rho}_n = \alpha_n + \beta_n + \gamma_n + \delta_n$. 
Let us first comment on the large sample properties of $\alpha_n$. 
Since $\alpha_n$ is the i.i.d. average of terms with 
expectation $\rho$ for all $P \in \mathcal{P}_0$ we have that 
$\alpha_n \to_P \rho$ for all $P \in \mathcal{P}_0$.
The term $\hat{\beta}_n$ is dealt with 
similarly as in the proof of Theorem \ref{thm:rhoasympnormal}. 
For fixed $k, \ell = 1, \dots, q$ we have that
	\begin{align*}
	|\beta_{k \ell, n}| & \leq 
	\frac{1}{n}
    \sum_{i=1}^n
    \left|
    \phi_k(\hat{U}_{1, i}) - \phi_k(U_{1, i})
    \right|
    \left|
    \phi_\ell(\hat{U}_{2, i}) - \phi_\ell(U_{2, i})
    \right|
    \\
    & \leq 
    \|
	\phi_k \circ \hat{F}^{(n)}_{X \mid Z} 
	- 
	\phi_k \circ F_{X \mid Z}
	\|_{\infty} \cdot
	\|
	\phi_\ell \circ \hat{F}^{(n)}_{Y \mid Z} 
	-
	\phi_\ell \circ F_{Y \mid Z}
	\|_{\infty} \in \mathcal{O}_P(g_P(n) h_P(n))
	\end{align*}
where we have used Lemma \ref{lem:mylittlelemma}. From 
Assumption \ref{assump:Fhatconsistency} we get that 
$\beta_{k \ell, n} \to_P 0$ for each $k, \ell = 1, \dots, q$, and so 
$\beta_n \to_P 0$ for all $P \in \mathcal{P}_0$. The terms $\gamma_n$ and $\delta_n$ are analyzed 
similarly, so we only look at $\gamma_n$. We see that for $k, \ell 
= 1, \dots, q$,
	\begin{align*}
	| \gamma_{k \ell, n} | 
	& \leq 
	\frac{1}{n}
	\sum_{i = 1}^n 
	| \phi_k ( U_{1, i})|
	| \phi_{\ell}(\hat{U}_{2, i}) - \phi_{\ell} (U_{2, i}) | 
	\\
	& \leq 
	\|\phi_k \|_{\infty} \cdot
	\|
	\phi_\ell \circ \hat{F}^{(n)}_{Y \mid Z} 
	- \phi_\ell \circ F_{Y \mid Z}
	\|_{\infty} \in \mathcal{O}_P(h_P(n))
	\end{align*}
where we have used that $\| \phi_k \|_\infty < \infty$ since 
$\phi_k : [0, 1] \to \R$ is continuous and $[0, 1]$ is compact. Here we have 
used Lemma \ref{lem:mylittlelemma}, and we conclude that $\gamma_{k \ell, n} 
\to_P 0$ due to Assumption \ref{assump:Fhatconsistency}, which 
shows that $\gamma_n \to_P 0$ for all $P \in \mathcal{P}_0$. Conclusively, 
we have $\hat{\rho}_n \to_P \rho$ for all $P \in \mathcal{P}_0$. 
\hfill $\square$

\subsection{Proof of Corollary \ref{cor:power}} 

Assume that $P \in \mathcal{A}_0$ such that $\rho_{k \ell} \neq 0$ for some 
$k, \ell = 1, \dots, q$. Then we have
	\begin{align*}
	T_n = \| \Sigma^{-1 / 2} \hat{\rho}_n \Sigma^{-1 / 2} \|_F^2 \stackrel{P}{\to}
	\| \Sigma^{-1 / 2} \rho \Sigma^{-1 / 2} \|_F^2 > 0
	\end{align*}
for all $P \in \mathcal{A}_0$ because $\rho \neq 0$. Here we have used Theorem \ref{thm:rhohat_to_rho}. 
Therefore we obtain that
	\begin{align*}
	n T_n = n \| \Sigma^{-1 / 2} \hat{\rho}_n \Sigma^{-1 / 2} \|_F^2 \stackrel{P}{\to} 
	\infty 
	\end{align*}
for all $P \in \mathcal{A}_0$. This means that 
	\begin{align*}
	P ( n T_n > c) \to 1
	\end{align*}
as $n \to \infty$ for all $c \in \R$. From this we obtain that 
	\begin{align*}
	\liminf_{n \to \infty} 
	E_P(\hat{\Psi}_n) = 
	\liminf_{n \to \infty} P(n T_n > z_{1 - \alpha}) = 1
	\end{align*}
for all $\alpha \in (0, 1)$ whenever $P \in \mathcal{A}_0$. \hfill $\square$

\subsection{Proof of Proposition \ref{prop:UUtoXYZ}}

Assume $(U_1, U_2) \indep Z$ and $U_1 \indep U_2$. Then it also holds that 
$U_1 \indep U_2 \mid Z$, which gives $(U_1, Z) \indep (U_2, Z) \mid Z$. 
More explicitly we have
  \begin{align*}
  (F_{X \mid Z}(X \mid Z), Z) \indep (F_{Y \mid Z}(Y \mid Z), Z) \mid Z.
  \end{align*}
Transforming with the conditional quantile functions gives
  \begin{align*}
  Q_{X \mid Z}(F_{X \mid Z}(X \mid Z) \mid Z)
  \indep 
  Q_{Y \mid Z}(F_{Y \mid Z}(X \mid Z) \mid Z) \mid Z.
  \end{align*}
Since we assume throughout the paper that the conditional distributions 
$X \mid Z=z$ and $Y \mid Z=z$ are continuous for each $z \in [0, 1]^d$ we get that 
$(X, Z) \indep (Y, Z) \mid Z$ 
which reduces to $X \indep Y \mid Z$.  \hfill $\square$

\subsection{Proof of Theorem \ref{thm:uniformresults}}

We start by showing (i). Again we consider the decomposition 
$\hat{\rho}_n = \alpha_n + \beta_n + \gamma_n + \delta_n$ introduced 
in the proof of Theorem \ref{thm:rhoasympnormal}.
By the stronger 
condition of Assumption \ref{assump:uniformFhatconsistency} we immediately have 
that 
$\sqrt{n} \beta_n\to_{\mathcal{P}_0} 0, \sqrt{n} \gamma_n\to_{\mathcal{P}_0} 0$ 
and $\sqrt{n} \delta_n \to_{\mathcal{P}_0} 0$ by following the same arguments as in 
the proof of Theorem \ref{thm:rhoasympnormal}. The fact that 
$\sqrt{n} \alpha_n$ converges uniformly in distribution to a $\mathcal{N}(0, \Sigma \otimes \Sigma)$-distribution 
over $\mathcal{H}_0$ follows 
from the fact that the distribution of 
$(U_{1, i}, U_{2, i})_{i=1}^n$ is unchanged whenever $P \in \mathcal{H}_0$. 
By Lemma \ref{lemma:slutsky} we have that 
  \begin{align*}
  \sqrt{n} \hat{\rho}_n = \sqrt{n} \alpha_n + \sqrt{n} \beta_n 
  + \sqrt{n} \gamma_n + \sqrt{n} \delta_n \Rightarrow_{\mathcal{H}_0} 
  \mathcal{N}(0, \Sigma \otimes \Sigma) 
  \end{align*}
which shows part (i) of the theorem. Next we turn to part (ii) of the theorem. 
Analogously to the proof of Theorem \ref{thm:rhohat_to_rho} we have that 
$\beta_n\to_{\mathcal{P}_0} 0, \gamma_n\to_{\mathcal{P}_0} 0$ and 
$\delta_n \to_{\mathcal{P}_0} 0$ under 
Assumption \ref{assump:uniformFhatconsistency}. Now 
consider writing 
	\begin{align*}
	\alpha_{k \ell, n} = \frac{1}{n} \sum_{i=1}^n A_i 
	\quad \text{where} \quad 
	A_i = \phi_k(U_{1, i}) \phi_{\ell}(U_{2, i})
	\end{align*}
for $k, \ell = 1, \dots, q$. Then $(A_i)_{i=1}^n$ are i.i.d. 
with $E_P(A_i) = \rho_{k \ell}$ and 
	\begin{align*}
	V_P(A_i) = E_P(\phi_k(U_{1, i})^2 \phi_\ell(U_{2, i})^2) - \rho^2
	\leq \| \phi_k \|_{\infty} \| \phi_{\ell} \|_{\infty} < \infty
	\end{align*}
for all $P \in \mathcal{P}_0$. Therefore, for given $\epsilon > 0$, we have
by Chebyshev's inequality that
	\begin{align*}
	\sup_{P \in \mathcal{P}_0} P(|\alpha_{k \ell, n} - \rho_{k \ell}| > \epsilon) 
  \leq 
  \sup_{P \in \mathcal{P}_0} \frac{V_P(\frac{1}{n}\sum_{i=1}^n A_i)}{\epsilon^2}
  =
  \sup_{P \in \mathcal{P}_0} \frac{V_P(A_i)}{n \epsilon^2}
  \leq \frac{||\phi_1||^2_{\infty} ||\phi_2||^2_\infty}{n \epsilon^2} \to 0
	\end{align*}
for $n \to \infty$ which shows that $\alpha_n \to_{\mathcal{P}_0} \rho$. 
From this we get $\hat{\rho}_n \to_{\mathcal{P}_0} \rho$ as wanted.

\subsection{Proof of Corollary \ref{cor:uniformlevel}}

Note that due to Theorem \ref{thm:uniformresults} (i) we have that 
$n T_n \Rightarrow_{\mathcal{H}_0} \chi^2_{q^2}$ under Assumption 
\ref{assump:uniformFhatconsistency} using the same argument as in 
the proof of Corollary \ref{cor:chisqteststat}. 
Then the result is obtained by the same argument as in the proof of 
Corollary $\ref{cor:level}$ by noting that 
$\sup_{P \in \mathcal{H}_0}|F_{n T_n}(t) - \Phi(t)| \to 0$ as 
$n \to \infty$ for all $t \in \R$. \hfill $\square$

\subsection{Proof of Corollary \ref{cor:uniformpower}} 

Let $\lambda > 0$ be fixed. By Theorem \ref{thm:uniformresults} (ii) 
we have $\hat{\rho}_n \to_{\mathcal{A}_\lambda} \rho$ where 
$|\rho_{k \ell}| > \lambda > 0$ for some $k, \ell = 1, \dots, q$. 
Therefore $\inf_{P \in \mathcal{A}_\lambda} |\rho_{k \ell}| \geq \lambda > 0$ 
and so 
	\begin{align*}
	T_n = \| \Sigma^{-1/2} \hat{\rho}_n \Sigma^{-1/2} \|_F^2 
	\to_{\mathcal{A}_\lambda} \| \Sigma^{-1/2} \rho \Sigma^{-1/2} \|_F^2  > 0 
	\end{align*}
since $\inf_{P \in \mathcal{A}_\lambda}|\rho^P_{k \ell}| > 0$ and $\Sigma^{-1/2}$ is positive definite. Therefore $n T_n \to_{\mathcal{A}_{\lambda}} \infty$, and 
so we have
	\begin{align*}
	\inf_{P \in \mathcal{A}_\lambda} 
	P(n T_n > c)
	& =
	\inf_{P \in \mathcal{A}_\lambda} 
	(1 - P(n T_n \leq c))
	\\
	& = - \sup_{P \in \mathcal{A}_\lambda} (P(n T_n \leq c)) - 1) 
	\to 1
	\end{align*}
as $n \to \infty$ for all $c \in \R$. From this we have
	\begin{align*}
	\liminf_{n \to \infty} \inf_{P \in \mathcal{A}_\lambda}
	E_P(\hat{\Psi}_n) = 
	\liminf_{n \to \infty} \inf_{P \in \mathcal{A}_\lambda}
	P(n T_n > z_{1 - \alpha}) = 1
	\end{align*}
for all $\alpha \in (0, 1)$. \hfill $\square$

\section{Modes of Stochastic Convergence} \label{appendix:modesofconvergence}

Let $\mathcal{M}$ denote some class of distributions. 
We start by defining the notions of small and big O in probability.

\subsection{Small and big-O in probability} \label{appendix:smallbigO}

All sequences $(a_n)$ and $(b_n)$ below are assumed to 
be non-zero.

\begin{definition}
Let $(X_n)$ and $(a_n)$ be sequences of random variables in $\R$.
If for every $\epsilon > 0$
  \begin{align*}
  \sup_{P \in \mathcal{M}} P(|X_n / a_n| > \epsilon) \to 0 
  \end{align*}
for $n \to \infty$ then we say that $X_n$ is small 
O of $a_n$ in probability uniformly over $\mathcal{M}$ and write 
$X_n \in o_\mathcal{M}(a_n)$. If for every $\epsilon > 0$ there is 
$M > 0$ such that
  \begin{align*}
  \sup_{n \in \mathbb{N}} \sup_{P \in \mathcal{M}} P(|X_n / a_n| > M) < \epsilon
  \end{align*}
then we say that $X_n$ is big 
O of $a_n$ in probability uniformly over $\mathcal{M}$ and write 
$X_n \in \mathcal{O}_{\mathcal{M}}(a_n)$.
\end{definition}

When $X_n \in \mathcal{O}_\mathcal{M}(a_n)$ we also say that $X_n$ 
is stochastically bounded by $a_n$ uniformly over $\mathcal{M}$.
When $X_n \in o_{\mathcal{M}}(1)$ we will typically write 
$X_n \to_\mathcal{M} 0$.

\begin{lemma} \label{lemma:bigO1}
Let $(X_n), (a_n)$ and $(b_n)$ be sequences of random 
variables in $\R$ such that
$X_n \in \mathcal{O}_{\mathcal{M}}(a_n)$. 
Then it holds that
$b_n X_n \in \mathcal{O}_{\mathcal{M}}(a_n b_n)$.
\end{lemma}

\begin{lemma} \label{lemma:bigO2}
Assume that $X_n \in \mathcal{O}_\mathcal{M}(a_n)$
and $Y_n \in \mathcal{O}_\mathcal{M}(b_n)$. 
Then $X_n Y_n \in \mathcal{O}_\mathcal{M}(a_n b_n)$. 
\end{lemma}

\begin{lemma} \label{lemma:bigO3}
Assume $X_n \in \mathcal{O}_\mathcal{M}(a_n)$ and that $a_n \in o(1)$. 
Then $X_n \in o_\mathcal{M}(1)$. 
\end{lemma}

\begin{lemma} \label{lemma:bigO4}
Assume that $X_n \in o_\mathcal{M}(1)$ and that $|X_n| \leq C$ for all $n \geq 1$ 
for a constant $C$ that does not depend on $P$. Then 
$\sup_{P \in \mathcal{M}}E_P|X_n| \to 0$ for $n \to \infty$.  
\end{lemma}

We now turn to uniform convergence in distribution.

\subsection{Uniform convergence in distribution} \label{appendix:uniformconvergencedistribution}

We follow \cite{kasy2019uniformity} and \cite{bengs2019uniform}.

\begin{definition}
Let $X, X_1, X_2, \dots$ be real valued random variables with distribution determined by 
$P \in \mathcal{M}$. If it holds that
  \begin{align*}
  \sup_{P \in \mathcal{M}}
  |
  E_P(f(X_n)) - E_P(f(X))
  | \to 0
  \end{align*}
for $n \to \infty$ for all functions $f : \R \to \R$ that are bounded and continuous, then we
say that $(X_n)$ converges uniformly in distribution to $X$ over $\mathcal{M}$. 
In this case we write $X_n \Rightarrow_{\mathcal{M}} X$. 
\end{definition}

\begin{lemma}[Uniform Slutsky's Lemma] \label{lemma:slutsky}
Assume that $X_n \Rightarrow_\mathcal{M} X$ and that $Y_n \to_\mathcal{M} 0$. Then 
$X_n + Y_n \Rightarrow_\mathcal{M} X$.
\end{lemma}

\begin{proof}
See \cite{bengs2019uniform} Theorem 6.3.
\end{proof}

\vskip 0.2in
\bibliography{biblio}

\end{document}